\documentclass[12pt]{amsart}
\usepackage{amsfonts}
\usepackage{latexsym,amsmath,amsthm,amssymb,amsfonts}
\usepackage{graphicx}
\numberwithin{equation}{section}

\def\S{\Sigma}
\def\CS{\mathcal S}

\def\endproof{$\hfill\Box$\\}

\def\s{\,\,\,\,}
\def\R{\mathbb{R}}
\def\P{\mathbb{P}}
\def\C{\mathbb{C}}

\numberwithin{equation}{section}
\newtheorem{theorem}{Theorem}[section]
\newtheorem{lem}[theorem]{Lemma}
\newtheorem{thm}[theorem]{Theorem}
\newtheorem{pro}[theorem]{Proposition}

\newtheorem{defi}[theorem]{Definition}

\newtheorem{exa}[theorem]{Example}
\def\dint{\displaystyle{\int}}

\begin{document}
\title[Homotopy classes of harmonic maps]{Homotopy classes of harmonic maps of the stratified 2-spheres and applications to geometric flows}

\author{Jingyi Chen and Yuxiang Li}
\address{ Department of Mathematics\\ The University of British Columbia, Vancouver, BC V6T1Z2, Canada}
\email{jychen@math.ubc.ca}
\address{Department of Mathematical Sciences, Tsinghua University, Beijing 100084, China}
\email{yxli@math.tsinghua.edu.cn}
\thanks{The first author acknowledges the partial support from NSERC and wishes to thank Tsinghua University where a part of this work was carried out during his visits. The second author is partially supported by NSFC. Both authors are grateful to Jiajun Wang for useful discussion on topology of 3-manifolds.}
\date{}
\maketitle
\begin{abstract}
We show that the set of harmonic maps from the 2-dimensional stratified spheres with uniformly bounded energies contains only finitely many homotopy classes. We apply this result to construct infinitely many harmonic map flows and mean curvature flows of 2-sphere in the connected sum of two closed 3-dimensional manifolds $M_1\not=S^3$ and $M_2\not=S^3,\R\P^3$, which must develop finite time singularity. 
\end{abstract}

\section{Introduction}

The fundamental existence result \cite{S-U} of Sacks-Uhlenbeck asserts that there exists a set of free homotopy classes of maps from $S^2$ to a
compact manifold $N$ so that the elements in these classes form a generating set for $\pi_2(N)$ acted on by $\pi_1(N)$ and each class contains a minimizing harmonic map $S^2\to N$. This can be used to show that there are only finitely many distinct free homotopy classes in
$\pi_2(N)$ which are representable by maps in $W^{1,2}(S^2,N)$ with energies stay below a given positive constant \cite{Sch-Wol}. For homotopy classes, one needs to count contribution of $\pi_1(N)$ by fixing a base point, and the set $\{ f\in W^{1,2}(S^2,N): E(f)<C\}$ where $C$ is a positive constant and $E(f)$ is the Dirichlet energy of $f$ may contain infinitely many homotopy classes. An example describing such behavior is first constructed in \cite{D-K} and then modified in \cite{Li-Wang} to the compact setting by mapping a portion of the domain surface to curves whose lengths go to infinity while the total energies remain uniformly bounded.


For harmonic maps, if the domain is a single $S^2$ and the energies are uniformly bounded, it follows from Parker's bubble tree convergence theorem \cite{Parker} that the limiting map (of a subsequence) is from a  stratified sphere, i.e. a chain of 2-spheres obtained by pinching finitely many loops on $S^2$ to points (see definition in Section 2), and consequently there are only finitely many homotopy classes for harmonic maps from $S^2$ under a fixed energy level. 

There are interesting geometric situations, e.g. harmonic map heat flow and mean curvature flow,  where we need to consider a sequence of harmonic maps from stratified spheres, possibly different, and study the homotopy  of these maps. 

In this paper we prove: 
\begin{thm}\label{str}
Let $N$  be a compact Riemannian manifold without boundary.
 Then for any $C>0$, the set of harmonic maps from stratified 2-spheres to $N$ with energies $\leq C$
contains only finitely many homotopy class.
\end{thm}

To keep track of the homotopy classes when a sequence of harmonic maps from (possibly different) stratified spheres converge as a bubble tree, we reformulate the construction of bubble trees \cite{Parker}, \cite{Parker-Wolfson} and introduce {\it a normalization} of a stratified sphere (see Definition \ref{homo}), which is, roughly, a surjective continuous map from $S^2$ to a stratified sphere $\Sigma$ and collapses finitely many multiply connected regions arising from the blowup points to the singular points of $\Sigma$. This procedure allows us to use the homotopy class of a normalization to capture the homotopy of a map from a stratified sphere.\\

There are interesting applications of Theorem \ref{str} to geometric flows. In this regard, 
it is useful to find manifolds which have infinitely many homotopy classes representable by mappings from $S^2$ with uniformly bounded energies. 
In Theorem \ref{3-manifold}, we show that if $M_1,M_2$ are closed 3-manifolds and $M_1\not=S^3,M_2\not=S^3,{\R\P^3}$ and $M=M_1\# M_2$ is the connected sum of $M_1,M_2$ then there are infinitely many smooth maps from $S^2$ to $M$ with uniformly bounded energies and are mutually nonhomotopic. The key topological result is Proposition \ref{h.sequence} derived in Appendix.

As applications of Theorem \ref{str} and Theorem \ref{3-manifold}, we prove finite time blowup for the harmonic map flow and the mean curvature flow of $S^2$ in the 3-manifold $M$ above. The formation of finite time singularities is due to the topology of $M$. For the harmonic map flow, these new examples and their construction are completely different from the one in \cite{C-D-Y} by Chang-Ding-Ye.

\begin{thm}\label{flow} 
Let   $M_1\neq S^3$ and $M_2\neq S^3, \R\P^3$ be closed 3-dimensional Riemannian manifolds and $M=M_1\# M_2$. 
Then 

(A) there exist infinitely many smooth maps $u_k:S^2\to M$ such that $u_i,u_j$ are not homotopic for any $i\neq j$, $\sup_k E(u_k)<\infty$ and 
the harmonic map flow that begins at $u_k$ develops a finite time singularity. 

(B) there exist infinitely many embeddings $w_k:S^2\to M$ such that $w_i,w_j$ are nonhomotopic for any $i\neq j$ with uniformly bounded area, and the mean curvature flow of $S^2$ initiated at $w_k$ develops a singularity in a finite time.
\end{thm}

To apply Theorem \ref{str} in the proof of Theorem \ref{flow}, for each of the two flows, we need a convergence result for a subsequence along the flow with the chosen initial data. 
It is shown by Qing-Tian in \cite{Qing-Tian} that for the bubbling of a Palais-Smale sequence for the Dirichlet energy functional with a uniform $L^2$-bound on the tension fields, there are no necks between bubbles in the limit. When the domain surface is $S^2$, the limit is a harmonic map from a stratified 2-dimensional sphere. Hence, absence of harmonic maps from any stratified spheres with energies less than a given bound $C$ in certain homotopy classes implies that the harmonic map flow starting from an initial map with energy below $C$ in such homotopy classes must develop finite time singularities, as pointed out already in Introduction of \cite{Qing-Tian}.  For the mean curvature flow, we need to invoke the compactness theorem 
in \cite{C-L} to extract a limiting harmonic stratified sphere if the flow exists for all time.


\section{Stratified Riemann surfaces and their normalizations} In this section, we introduce a notion, called normalization,  that provides a systematic way to associate a mapping  $S^2\to \Sigma$, for any given stratified sphere $\Sigma$. The normalization allows us to study homotopy classes of mappings from (possibly different) stratified spheres to a manifold $N$ by using mappings from a single $S^2$ to $N$.


Let $(M, d)$ be a connected compact metric space. As in \cite{C-T}, we say $M$ is
a  {\em stratified Riemann surface with singular set $P$}
if $P\subset M$ is finite set such that
\begin{enumerate}
\item  $(M\backslash P,d)$ is a smooth Riemann
surface without boundary (possibly disconnected) and $d$ is a smooth
metric $h=d|_{M\backslash P}$;
\item  For each $p\in P$, there is $\delta$, such that
$B_\delta(p)\cap P=\{p\}$ and $B_\delta(p)\backslash\{p\}
=\cup_{i=1}^{m(p)}\Omega_i$,
where $1<m(p)<+\infty$, and each $\Omega_i$ is topologically
a disk with its center deleted. Moreover,
$h$ can be extended on each $\Omega_i$ to be a smooth metric on the disk.
\end{enumerate}
A connected component of $M\backslash P$ is called a {\it component of $M$}. The genus of $M$ is
$$
g(M)={1-\frac{1}{2}\chi(M)+\frac{1}{2}\sum\limits_{p\in P}(m(p)-1)}.
$$
When $g(M)=0$, $M$ is called a {\it stratified sphere}.
A stratified sphere can be obtained by shrink finitely many disjoint embedded loops in $S^2$ to points.

A continuous map $u$ from a stratified surface $M$ to $N$  {\it harmonic}
if $u$ is harmonic  on each component of $M$. The map $u$ may be trivial on some components. 

\subsection{Normalization of stratified spheres} A domain in a smooth sphere is a {\it multiply connected collar}
if it is topologically the sphere minus $k$  open disks $D_1$, $\cdots$, $D_k$
for some positive integer $k$, where $D_i\cap D_j=\emptyset$ for any $i\neq j$. We will use the multiply connected collars to describe 
the neck regions in blowup analysis. 

\begin{defi}\em{Let $\S$ be a stratified sphere with $m$ singular points
$P=\{p_1, \cdots, p_m\}$. A continuous surjective map
$\phi: S^2\rightarrow \S$ is called {\it a normalization of $\S$}, if there exist
multiply connected collars  $\Omega_1$, $\cdots$,
$\Omega_m\subset S^2$, such that
\begin{enumerate}
\item  $\Omega_i\cap \Omega_j=\emptyset$, $i\not= j$;
\item $\phi$ is injective on $S^2\backslash\bigcup\limits_{i=1}^m\Omega_i$
and $\phi(S^2\backslash\bigcup\limits_{i=1}^m\Omega_i)\bigcap P=\emptyset$;
\item $\phi(\Omega_i)=p_i$.
\end{enumerate}
We denote a normalization by $\phi:(S^2; \Omega_1,\cdots,\Omega_m)\to (\S;p_1,\cdots,p_m)$. }
\end{defi}
When $m=0$, i.e. $P=\emptyset$, $\S$ is a smooth $S^2$ and any homeomorphism from $S^2$ to $\S$ is a normalization.

\begin{figure}[!ht]
\includegraphics{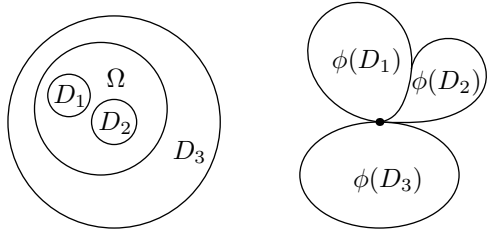}
\caption{An example of normalization: $\Omega=S^2\backslash \cup^3_{k=1}D_k$ is a multiply connected collar, $\S=\phi(D_1\cup D_2\cup D_3)\cup \{p\}$.}
\end{figure}

The following result is useful to construct a normalization for a given stratified sphere. 
\begin{lem}\label{nor}
Let $\Omega_1$, $\cdots$, $\Omega_m$ be multiply connected collars in $S^2$ with $\Omega_i\cap\Omega_j
=\emptyset$ for $i\neq j$. Let $f$ be a continuous map from $S^2$ to a stratified
sphere $\S$ with singular point set $P$. If $f(\Omega_i)\in P$ for each $i=1,\cdots,m$ and $f$ is bijective
from $S^2\backslash\cup_{i=1}^m\Omega_{i}$ to $\S\backslash P$,
then $P$ has $m$ elements and $f:(S^2;\Omega_1,\cdots,\Omega_m)\rightarrow (\S; f(\Omega_1),\cdots,f(\Omega_m))$ is a normalization.
\end{lem}

\proof We argue by induction on the number of components of $\S$.

When $\S$ has only one component, $P$ is empty and the lemma obviously holds.
Suppose the result is valid for any $\S$ with $k-1$ components.

Now assume $S$ has $k$ components. Since $\S$ is a stratified sphere, it has a component $S_1$  whose topological closure $\overline{S_1}$  intersects $P$ only at one point $p$.  By assumption, $f^{-1}$ is bijective from the topological disk $\overline{S_1}\backslash \{p\}=S_1$ to its image $U=f^{-1}(S_1)$, so $U$ is a topological disk.
Without loss of generality, we assume $\Omega_m\subseteq f^{-1}(p)$.
Since $f$ is bijective from  $S^2\backslash \cup \,\Omega_i$ to $\S\backslash P$, we have $U\subseteq S^2\backslash \cup\,\Omega_i$ and $\partial U\subseteq\partial\Omega_m$. Set $\S'=(\S\backslash S_1)\cup\{p\}$. Clearly, $\S'$ is a stratified sphere with $m-1$ components and let $P'$ be its singular set.

If $p$ is not in  $P'$, then $\Omega_m$ is an annulus,
and $\Omega_m\cup U$ is a disk. Since $U\cup\Omega_m$ does not intersect $\Omega_i$ for any $i<m$, we can take a disk type neighbourhood $V$ of $U\cup\Omega_m$ in $S^2$, such that $V$
does not intersect $\Omega_i$ for any $i<m$. Then we can extend
$f|_{S^2\backslash V}$ to a map $f':S^2\rightarrow \S'$ such that
$f'|_{S^2\backslash V}=f$ and $f'$ is bijective from $V$ onto the disk $f(V\backslash U)\cup\{p\}$.
Then $f'$ is bijective from $S^2\backslash \cup_{i=1}^{m-1}\Omega_i$ to $\S'\backslash P'$. By the induction hypothesis,  $\S'$ has exact
$m-1$ singular points $p_1$, $\cdots$, $p_{m-1}$, and we may assume
$f'(\Omega_i)=f(\Omega_i)=p_i$ for $i<m$. Hence
$f:(S^2;\Omega_1,\cdots,\Omega_m)\rightarrow (\S; p_1,\cdots,p_{m-1},p)$
is a normalization of $\S$.

If $p$ is a singular point of $\S'$, then
$\Omega_m = D\backslash \cup^{m'}_{i=1} D_{i}$ with $m'>1$, where $D, D_1,\cdots,D_{m'}$ are disks in $S^2$ and $D_i\cap D_j=\emptyset$ if $i\not=j$ (equivalently, $\Omega_m$ is $S^2$ with at least 3 non-intersecting disks removed). Note $U$ is one of the $D_i$'s. 
Define $\Omega_m'=\Omega_m\cup U$ and
$$f'(x)=\left\{\begin{array}{ll}
                f(x)&x\notin U\\
               p&x\in U.\end{array}\right.$$
Then $f'$ is bijective from $S^2\backslash (\cup_{i=1}^{m-1}\Omega_i \cup\Omega_m')$
to $\S'\backslash P'$. By the induction hypothesis,  $\S'$ has exact
$m-1$ singular points $p_1$, $\cdots$, $p_{m-1}$ except $p$, and we may assume
$f'(\Omega_i)=f(\Omega_i)=p_i$. Hence
$f:(S^2;\Omega_1,\cdots,\Omega_m)\rightarrow (\S;p_1,\cdots,p_{m-1},p)$
is a normalization of $\S$. Now the induction is complete.
\endproof

For a given stratified sphere, its normalization is unique up to homeomorphisms:
\begin{lem}\label{Lemma}
Let $\S$ be a stratified sphere with $P=\{p_1,\cdots,p_m\}$ and let
$$\phi:(S^2;\Omega_1,\cdots,\Omega_m)\rightarrow
(\S;p_1,\cdots,p_m) \,\,
\mbox{and} \,\,\phi':(S^2;\Omega_1',\cdots,\Omega_m')\rightarrow
(\S;p_1,\cdots,p_m)$$
be two
normalizations of $\S$. Then there is a homeomorphism $\psi:S^2\rightarrow S^2$
such that $\phi'=\phi\circ\psi$.
\end{lem}

\proof

We will argue by induction on the number $m$. The lemma obviously holds when $m=0$ as both $\phi,\phi'$ are homeomorphisms from $S^2$ to $S^2$, and assume it also holds for all $m\leq k-1$.

For $m=k$,  take $p\in P$. Assume  $\S\backslash\{p\}$ has $n$ connected components
$U_1$, $\cdots$, $U_{n}$ (each $U_i$ can be a union of 2-spheres). Since $\S$ has genus 0, there is no closed chain of $S^2$'s in $\S$, therefore  $\S_i=U_i\cup\{p\}$ is a stratified sphere with no more than $k-1$ points in its singular set $P_i$, by noting that $p$ is not a singular point of $\S_i$, for each $i=1,\cdots, n$. Then $\phi^{-1}(p)$ is a multiply connected collar with $n$ boundary circles, i.e.,
$S^2\backslash \phi^{-1}(p)$ has $n$ connected components $D_1,\cdots,D_n$ which are all topological disks. Note that $\phi(\overline{D_i})= \S_i$ and $\phi(\partial D_i)=p$. Let $\alpha_i$  be a continuous map from
$\overline{D_i}$ to $S^2$ such that $\alpha_i(\partial D_i)=N$ and $\alpha_i$ is bijective
from $D_i$ to $S^2\backslash \{N\}$, where $N$ is the north pole of $S^2$. Define $\phi_i:S^2\to \S_i$ by 
$$\phi_i(x)=\left\{\begin{array}{ll}
          \phi\circ(\alpha_i)^{-1}(x),&x\in S^2\backslash N\\
          p,&x=N,
\end{array}\right.$$
then it is easy to check that $\phi_i$ is a normalization of $\S_i$ corresponding to the multiply connected collars $\alpha(\Omega_i)$ subject to 
$\phi(\Omega_i)=p_i \in P_i$.

By the same argument,  $\phi'^{-1}(p)$ is a multiply connected collar with $n$ boundary curves, and
 we  define  a normalization $\phi'_i$ of $S_i$ as above. Then, by the induction hypothesis,
we can find a homeomorphism $\psi_i:S^2\to S^2$ such that $\phi_i'=\phi_i\circ\psi_i$, for each $i=1,\cdots,n$.

Since both $S^2\backslash \phi^{-1}(p)$ and $S^2\backslash \phi'^{-1}(p)$  are unions of $n$ disjoint topological disks
$D_1,\cdots,D_n$ and $D'_1,\cdots,D'_n$, respectively,
 we can find  a homeomorphism $\psi :S^2\rightarrow S^2$, such that
(1) $\psi(\phi'^{-1}(p))=\phi^{-1}(p)$ and (2) $\psi(D_i')=D_{i}$ and $\psi|_{D_i'} = \alpha_i^{-1}\circ\psi_i\circ\alpha_i'$.
We now show $\phi\circ\psi = \phi'$:

(A) For any $x\in \phi'^{-1}(p)$,  $\phi'(x)=p$ and $\psi(x)\in\phi^{-1}(p)$ by (1), hence  $\phi\circ\psi(x)=\phi'(x)$.

(B) For any $x\in D_i'$, we have $x\not\in\phi'^{-1}(p)$.
If $\psi_i(\alpha'_i(x))=N$ then  by (2)
$$
\phi\circ\psi(x)= \phi\circ\alpha_i^{-1}\circ\psi_i\circ\alpha'_i(x)=\phi\circ\alpha^{-1}_i(N)
=\phi(\partial D_i)=p
$$
So $\psi(x)\in\phi^{-1}(p)$ and in turn $x\in\phi'^{-1}(p)$ by (1). This is contradiction. Therefore
$$
\phi\circ\psi(x)=\phi\circ\alpha_i^{-1}\circ\psi_i\circ\alpha'_i(x)
=\phi_i\circ \psi_i\circ\alpha'_i(x)=\phi'_i\circ\alpha'_i(x)=\phi'(x)
$$
provided $\alpha_i'(x)\not=N$. But $\alpha'_i(x)=N$ would imply $x\in\partial D'_i$ hence $x\in \phi'^{-1}(p)$ which is impossible. In conclusion, we have 
$\phi\circ\psi=\phi'$. 
\endproof

\subsection{Homotopy classes of mappings from stratified spheres} 
We generalize the standard notion of homotopy equivalence of two maps from $S^2$ to $N$ to two maps from stratified spheres, possibly different, to $N$.

\begin{defi}\label{homo}{\em
Let $\S_1$ and $\S_2$ be two stratified spheres. Let $u_i$ be continuous maps from $\S_i$ to $N$, $i=1,2$.
We say $u_1$ is {\it homotopic} to $u_2$ if there exist normalizations $\phi_1,\phi_2$ of $S_1,S_2$, respectively,
such that $u_1\circ\phi_1$ is homotopic to $u_2\circ\phi_2$ as maps from $S^2$ into $N$. We write $u_1\sim u_2$  if $u_1,u_2$ are homotopic.}
\end{defi}
 Two maps being homotopic is independent of the choice of normalizations:  Suppose $u_1\sim u_2$. Then there exist normalizations $\phi_i$ of $\S_i$ so that $u_1\circ\phi_1\sim u_2\circ\phi_2$. Let $\phi'_1,\phi'_2$ be two arbitrary normalizations of $\S_1,\S_2$, respectively. By Lemma \ref{Lemma}, there exist homeomorphisms $\psi_i:S^2\to S^2$ such that
$\phi_i=\phi_i'\circ\psi_i$ for $i=1,2$. Then  $u_1\circ \phi_1'\circ \psi_1\sim u_2\circ\phi_2'\circ\psi_2$. Since any homeomorphism of $S^2$ to itself is homotopic to the identity map, it follows $u_1\circ \phi_1'\sim u_2\circ\phi_2'$.

\begin{lem}\label{homotopy}
Suppose $\S$ is a stratified sphere with components $S_1$, $\cdots$, $S_m$ and each $S_i$ is a 2-sphere $S^2$. Let  $\hat{\S}=\cup_{i=1}^m \hat{\S}_i$ be a stratified sphere, where each $\hat{\S}_i$
is a stratified sphere and $\hat{\S}_i\cap \hat{\S}_j$ consists of singular points of $\hat{\S}$ for any $i\not= j$. Suppose
that $f:\S\rightarrow \hat{\S}$ is a continuous map such that $f|_{S_i}$ is a normalization of $\hat{\S}_i$ and $f$ maps each
singular point of  $\S$ to a singular point of $\hat{\S}$. Then  $u\sim u\circ f$ for any
$u\in C^0(\hat{\S},N)$.
\end{lem}
\proof It suffices to show that there is a normalization $F:S^2\rightarrow \S$
such that $f\circ F$ is a normalization of $\hat{\S}$, because $(u\circ f)\circ F=u\circ(f\circ F)$.
We argue by the induction on $m$. When $m=1$, $\S=S^2$ and $f:\S\to \hat{\S}$ is a normalization of $\hat{\S}$, we may take $F$ to be the identity map of $S^2$. Assume the claim is true for $m=k$.
For $m=k+1$, without loss of generality, we assume
$S_1$ is a component of $\S$ which contains only one singular point $p$ of $\S$.
Set $\S'=\cup_{i=2}^{k+1}S_i$ and $\hat{\S}'=\cup_{i=2}^{k+1}
\hat{\S}_i$.  Let  $f'=f|_{\S'}$. Then $f':\S'\rightarrow \hat{\S}'$ satisfies the conditions
in the lemma. Thus, by the induction hypothesis,  we can find a normalization $F'$ of $\S'$ such that
$f'\circ F'(S^2,\Omega_1,\cdots,\Omega_\alpha)\rightarrow (\hat{\S}',p_1,\cdots,p_\alpha)$
is a normalization.

Next, we construct $F$ on $\S$.  Note that $f|_{S_1}$ is a normalization of $\hat{\S_1}$ by assumption. Set
$$f_1=f|_{S_1}:(S_1;U_1,\cdots,U_\beta)\rightarrow (\hat{\S}_1;q_1,\cdots,q_\beta).$$
We have four cases:

Case 1: $p$ is not a singular point of $\S'$ and $f_1(p)$ is not a  singular point of $\hat{\S}_1$.

In this case, there exists a unique $x_0\in S^2=S_1$ with $F'(x_0)=p$.
Pick $r>0$ such  that $F'$ is bijective on $B_{4r}(x_0)$. Let $\varphi:B_{3r}(x_0)\backslash B_{2r}(x_0)\to B_{3r}(x_0)$ be a continuous map such that
\begin{enumerate}
\item $\varphi(\partial B_{2r}(x_0))=\{x_0\}$,
\item $\varphi$ is bijective from
$B_{3r}(x_0)\backslash \overline{B_{2r}(x_0)}$ to $B_{3r}(x_0)\backslash\{x_0\}$,
\item $\varphi$ is the identity map on $\partial B_{3r}(x_0)$.
\end{enumerate}
Let $\phi:\overline{B_r(x_0)}\to S_1$ be a continuous map such that
$\phi(\partial B_r(x_0))=\{p\}$ and $\phi$ is bijective from $B_r(x_0)$
onto $S_1\backslash\{p\}$. Set
$$
F=\left\{\begin{array}{ll}
                 F'(x)&x\in S^2\backslash B_{3r}(x_0)\\
                F'\circ\varphi(x)&x\in B_{3r}\backslash B_{2r}(x_0)\\
               p&x\in B_{2r}\backslash B_{r}(x_0)\\
               \phi(x)&x\in B_r(x_0).
\end{array}\right.
$$
Obviously, $F$ is a normalization of $\S$. It follows from Lemma \ref{nor} that
\begin{equation*}
\begin{split}
f\circ F: \s &(S^2;\Omega_1,\cdots,\Omega_\alpha,B_{2r}\backslash B_r(x_0),\phi^{-1}(U_1),\cdots,\phi^{-1}(U_{\beta}))\\ 
&\rightarrow(\hat{\S};p_1,\cdots,p_\alpha,f(p),q_1,\cdots,q_{\beta})
\end{split}
\end{equation*}
is a normalization of $\hat{\S}$. \\

Case 2: $p$ is not a singular point of $\S'$ but $f_1(p)$ is  a  singular point of $\hat{\S}_1$.

In this case, there exists a unique $x_0\in S^2$ with $F'(x_0)=p$. However,
we have $f_1(p)\in\{q_1,\cdots,q_\beta\}$.
We choose $x_0$, $r$, $\varphi$, $\phi$ and $F$ as before.
We may assume $f_1(p)=q_1$ and modify $f_1$ such that
$p\notin \partial U_1$.  Set $\Omega_{\alpha+1}=(\overline{B_{2r}(x_0)}\backslash B_r(x_0))
\cup \phi^{-1}(U_1)$. Note $\phi^{-1}(U_1)\subset \overline{B_r(x_0)}$. Then $\Omega_{\alpha+1}$ is a multiply connected collar homeomorphic to $U_1$. Clearly, $F$ is a normalization of $\S$ and by Lemma \ref{nor},
\begin{equation*}
\begin{split}
f\circ F:\s& (S^2;\Omega_1,\cdots,\Omega_\alpha,\Omega_{\alpha+1},
\phi^{-1}(U_2),\cdots,\phi^{-1}(U_{\beta})) \\
&\rightarrow
(\hat{\S};p_1,\cdots,p_\alpha,q_1,\cdots,q_{\beta})
\end{split}
\end{equation*}
is a normalization of $\hat{\S}$. \\

Case 3: $p$ is a singular point of $\S'$ and $f_1(p)$ is not a  singular point of $\hat{\S}_1$.

In this case, we may assume $f_1(p)=p_\alpha$ and $\Omega_\alpha=(f'\circ F')^{-1}(f'(p))$. Let $\Omega=F'^{-1}(p)$. It is easy to see $\Omega\subset\Omega_\alpha$: $\forall x\in \Omega$, $F'(x)=p$, then $f'(F'(x))=f'(p)$, so $x\in\Omega_\alpha$. Choose $r$ and $x_0$,  such  that  $B_{2r}(x_0)\subset \Omega$. Let $\phi$
be a continuous map from $\overline{B_r(x_0)}$ into $S_1$ such that
 $\phi(\partial B_r(x_0))=\{p\}$ and $\phi$ is bijective from
$B_r(x_0)$ to $S_1\backslash\{p\}$. Set
$$F=\left\{\begin{array}{ll}
                 F'(x)&x\in S^2\backslash B_r(x_0)\\
                 \phi(x)&x\in B_r(x_0)
                 \end{array}\right.$$
Then $F$ is a normalization of $\S$, and
\begin{equation*}\begin{split}
f\circ F:\s &(S^2;\Omega_1,\cdots,\Omega_{\alpha-1},\Omega_\alpha\backslash B_r(x_0),
\phi^{-1}(U_1),\cdots,\phi^{-1}(U_{\beta})) \\
&\rightarrow
(\hat{\S};p_1,\cdots,p_{\alpha-1},f(p),q_1,\cdots,q_{\beta})
\end{split}
\end{equation*}
is a normalization of $\hat{\S}$. \\

Case 4: $p$ is a singular point of $\S'$ and $f_1(p)$ is  a  singular point of $\hat{\S}_1$.

Without loss of generality, we may assume $f(p)=p_\alpha=q_1$, $\Omega_\alpha=(f'\circ \
F')^{-1}(f'(p))$.
Let $\Omega=F'^{-1}(p)\subset\Omega_\alpha$.
Choose $r$ and $x_0$,  such  that  $B_{2r}(x_0)\subset \Omega$.
Modify $f_1$ such that $p\notin \partial U_1$ and set $\Omega_\alpha'=
(\Omega_\alpha\backslash B_r(x_0))\cup\phi^{-1}(U_1)$.
Then $\Omega_{\alpha}'$ is a multiply connected collar. Define $F$ as in Case 3.
Then
\begin{equation*}\begin{split}
f\circ F:\s &(S^2;\Omega_1,\cdots,\Omega_{\alpha-1},\Omega_{\alpha}',
\phi^{-1}(U_2),\cdots,\phi^{-1}(U_{\beta})) \\
& \rightarrow
(\hat{\S};p_1,\cdots,p_\alpha,q_2,\cdots,q_{\beta})
\end{split}
\end{equation*}
is a normalization of $\hat{\S}$. \endproof

\section{Bubble tree construction }
This section is divided into three parts. Firstly, we assign an ordering for the blowup sequences according to their blowup rates and describe the adjacent ones. Secondly, we construct the neck region and show it is a union of disjoint multiply connected collars. This construction is convenient for using the normalizations to keep track the homotopy classes. Lastly, we collect known results which are needed for the bubble tree convergence.

\subsection{\bf Hierarchy of blowup sequences} We begin by recalling the following important results in \cite{S-U}.

\begin{pro}\label{epsilon} ($\epsilon$-regularity) There exists an $\epsilon_0>0$, such that for any harmonic map $u$ from the unit disk $D$ in ${\mathbb R}^2$ into $N$, if $\int_{D}|\nabla
u|^2\leq \epsilon_0^2$, then
$$\|\nabla^ku\|_{L^\infty(D_{1/2})}
\leq C(\epsilon_0,k,N)\|\nabla u\|_{L^2(D)},\s\forall k>0.$$
\end{pro}

\begin{pro}\label{gap}(Gap constant)
There exists $\tau>0$ which only depends on $(N,h)$, such that
there is no non-trivial harmonic map from $S^2$ into $N$, whose
energy is in $(0,\tau)$.
\end{pro}



Let $u_k$ be a sequence of harmonic maps from a compact surface $\Sigma$ with a Riemannian metric $g$ to a compact Riemannian manifold $N$ with $E(u_k)<\Lambda<\infty$.
We may assume $u_k\rightharpoonup u_0$ in $W^{1,2}
(\Sigma,g,N)$. The blowup set of the sequence $\{u_k\}$ is
$$\mathcal{S}=\mathcal{S}(\{u_k\})=\left\{p\in\Sigma:
\lim_{r\rightarrow0}\varliminf_{k\rightarrow+\infty}
\int_{B_r(p)}|\nabla_g u_k|^2\geq {\epsilon_0^2}\right\}.$$
Then we can find a subsequence, still denote by $\{u_k\}$, such that
for any $p\notin\mathcal{S}$, there is  $\delta>0$ with
$
\int_{B_{\delta}(p)}|\nabla u_k|^2\leq {\epsilon_0^2}.
$
Thus, we may assume that $\|\nabla^mu_k\|_{L^\infty
(\Omega)}<C(m,\Omega)$ for any $\Omega\subset\subset
\Sigma\backslash\mathcal{S}$, and that
 $u_k$ (possibly a subsequence) converges smoothly to $u_0$
in $\Omega$, as $N$ is compact.
The blowup set ${\mathcal S}$ consists of at most finitely many points which are called the (energy) concentration points of $\{u_k\}$. Near a point $p\in\mathcal{S}$,  we take an isothermal coordinate system $(D,x)$ centered at $p$ which is the only concentration point in $D$. The maps
$u_k$ can be regarded as harmonic maps from $D$ into $N$.

\begin{defi}\emph{
A sequence $\{(x_k,r_k): x_k\in D, r_k>0\}$ is a {\it blowup sequence of $\{u_k\}$ at $p\in{\mathcal S}$} if
$x_k\rightarrow 0$ and $r_k\rightarrow 0$, $v_k(x)=u_k(x_k+r_kx) \rightharpoonup v\,\, \mbox{in}\,\, W^{1,2}_{loc}(\R^2,N)$,
and
$v_k\rightarrow v \mbox{ in } C^\infty_{loc}(\R^2\backslash S(p))$ for some finite set $S(p)$, and $E(v)>0$.  }
\end{defi}
By removability  of singularity, $v$ extends to a harmonic map from $S^2$
into $N$. We call the non-constant harmonic map $v:S^2\to N$ a {\it bubble}. 

\begin{defi}\label{essen} \emph{
Two blowup sequences $\{(x_k,r_k)\}$ and $\{(x_k',r_k')\}$ of $\{u_k\}$ at $p\in{\mathcal S}$ are said to be
{\it essentially different} if one of the following happens
\begin{equation}\label{bubble1}
\frac{r_k}{r_k'}\rightarrow+\infty,\mbox{ or }
\frac{r_k'}{r_k}\rightarrow+\infty,
\s \mbox{or} \s\frac{|x_k-x_k'|}{r_k+r_k'}
\rightarrow+\infty.
\end{equation}
Otherwise, they are called {\it essentially same}. }
\end{defi}
In the sequel, we will write $(x_k,r_k)$ for a blowup sequence for simplicity. 
\begin{lem}\label{same}
If two blowup sequences $(x_k,r_k)$ and $(x_k',r_k')$ of $\{u_k\}$ at $p$ are essentially
same, then the bubbles $v,v'$ are the same in the sense $v=v'\circ L$, where $L:{\mathbb R}^2\to{\mathbb R}^2$ is a linear transformation.
\end{lem}
\proof Since (\ref{bubble1}) does not hold, after passing to a subsequence $k_n$ if necessary, we may assume
$$
\frac{r_{k_n}}{r_{k_n}'}\rightarrow \lambda\in(0,+\infty)\s \mbox{and}\s
\frac{x_{k_n}-x_{k_n}'}{r_{k_n}'}\rightarrow x_0.
$$
By the definition of blowup sequences, we have $v_k(x)=u_k(x_k+r_kx)\rightharpoonup v(x)$ and $v_k'(x)=u_k(x'_k+r'_kx)\rightharpoonup v'(x)$. Observe
$$
u_{k_n}(x_{k_n}+r_{k_n}x)=u_{k_n}\left(x_{k_n}'+r_{k_n}'\left(\frac{r_{k_n}}{r_{k_n}'}x+
\frac{x_{k_n}-x_{k_n}'}{r_{k_n}'}\right)\right).
$$
Letting $n\to\infty$ we then have
$$
v(x)= v'(\lambda x+x_0).
$$
We can then take $L(x)=\lambda x+x_0$.\endproof



\begin{lem}\label{different}
Let $\{u_k\}$ be a sequence of harmonic maps $D$ into $N$.
If $(x_k^\alpha,r_k^\alpha)$ are mutually essentially different  blowup sequences
 of $\{u_k\}$ at $0$,  where $\alpha=1,\cdots,m$,
then
$$
\varlimsup_{k\rightarrow+\infty}E(u_k,D)\geq m\tau,
$$
where $\tau$ is the gap constant in Proposition \ref{gap}.
\end{lem}
\proof Let $v_k^\alpha(x)=u_k(x_k^\alpha+r_k^\alpha x)$. We assume
$v_k^\alpha$ converges to $v^\alpha$ in $C_{loc}^\infty(\R^2\backslash \CS^\alpha)$,
where $\CS^\alpha$ is the set of concentration points of $\{v_k^\alpha\}$, which is finite. Set 
$$
U^\alpha_r=D_\frac{1}{r}\backslash \bigcup_{p\in \CS^\alpha}D_r(p),\s \Omega^\alpha_{r,k}=x_k^\alpha+r_k^\alpha U^\alpha_r.
$$
 As $v^\alpha$ is nontrivial, 
for any $\epsilon$, we can find $r$, such that
$$\lim_{k\rightarrow+\infty}E(v_k^\alpha,U_r^\alpha)=E(v^\alpha,U_r^\alpha)\geq \tau-\epsilon.$$
Then
$$
\lim_{k\rightarrow+\infty}E(u_k,\Omega^\alpha_{r,k})=
\lim_{k\rightarrow+\infty}E(v_k^\alpha,U^\alpha_{r})\geq \tau-\epsilon.
$$

By Definition \ref{essen}, we may assume one of the following alternatives holds:
\begin{enumerate}
\item $\frac{|x_k^\alpha-x_k^\beta|}{r_k^\beta+r_k^\alpha}\rightarrow+\infty$, which is equivalent to
$\frac{|x_k^\alpha-x_k^\beta|}{r_k^\beta}\rightarrow+\infty$ and
$\frac{|x_k^\beta-x_k^\alpha|}{r_k^\alpha}\rightarrow+\infty$,
\item  $\frac{r_k^\alpha}{r_k^\beta}\rightarrow 0$ and $\frac{x_{k_n}^\alpha-x_{k_n}^\beta}{r_{k_n}^\beta}\rightarrow x_{\alpha\beta}$  ($k_n$ is a sequence of positive integers, $n\in{\mathbb Z}^+$).
\end{enumerate}

When (1) holds, we have $D_{r_k^\alpha/r}(x_k^\alpha)
\cap D_{r_k^\beta/r}(x_k^\beta)=\emptyset$ when $k$ is sufficiently large.
Since $\Omega_{r,k}^\alpha\subset D_{r_k^\alpha/r}(x_k^\alpha)$ and
$\Omega_{r,k}^\beta\subset D_{r_k^\beta/r}(x_k^\beta)$, we conclude
$\Omega^\alpha_{r,k}\cap\Omega^\beta_{r,k}=\emptyset$ for large $k$.

When (2) holds, for any $\delta>0$, there exists $k(\delta)$, for all $k_n\geq k(\delta)$ we have
$$
E\left(v_{k_n}^\beta,D_{2\delta}(x_{\alpha\beta})\right)\geq E\left(v_{k_n}^\beta,D_\delta\left(\frac{x_{k_n}^\alpha-x_{k_n}^\beta}{r_{k_n}^\beta}\right)\right)
=E\left(u_{k_n},D_{r_{k_n}^\beta\delta}(x_{k_n}^\alpha)\right).
$$
Then, by noting that $\frac{r_k^\beta\delta}{r_k^\alpha}\rightarrow+\infty$ and writing
$u_{k_n}(x)=v^\alpha_{k_n}\left(\frac{x-x^\alpha_{k_n}}{r^\alpha_{k_n}}\right)$, we have
$$
\varliminf_{n\rightarrow+\infty}E\left(v_{k_n}^\beta,D_{2\delta}(x_{\alpha\beta})\right)\geq
\varliminf_{n\rightarrow+\infty}E\left(u_{k_n},D_{r_{k_n}^\beta\delta}(x_{k_n}^\alpha)\right)
\geq \lim_{n\rightarrow+\infty}E\left(v_{k_n}^\alpha,U_r^\alpha\right)\geq\tau-\epsilon.
$$
Thus, $x_{\alpha\beta}\in \CS^\beta$ and in turn $D_r(x_{\alpha\beta})\cap U_{r}^\beta=\emptyset$.
Then $D_{\frac{r}{2}}\left(\frac{x_{k_n}^\alpha-x_{k_n}^\beta}{r_{k_n}^\beta}\right)\cap U_{r}^\beta=\emptyset$ for large $k_n$, hence $D_{rr_{k_n}^\beta/2}(x_{k_n}^\alpha)\cap\Omega_{r,k_n}^\beta=\emptyset$ for large $k_n$.
Now it follows from
$$
\Omega_{k_n,r}^\alpha\subset D_{\frac{r_k^\alpha}{r}}(x_{k_n}^\alpha)\subset
D_{\frac{rr_{k_n}^\beta}{2}}(x_{k_n}^\alpha)
$$
 that
$
\Omega^\alpha_{r,k_n}\cap \Omega^\beta_{r,k_n}=\emptyset
$
for large $k_n$.

Therefore, we have established: for any $\alpha\not=\beta$, there exists an increasing sequence $\{k_n(\alpha,\beta):n\in{\mathbb N}\}$ of positive integers such that
\begin{equation}\label{differentbubble}
\Omega^\alpha_{r,k_n(\alpha,\beta)}
\cap\Omega^\beta_{r,k_n(\alpha,\beta)}=\emptyset.
\end{equation}
 Applying (\ref{differentbubble}) and the diagonal sequence procedure, we can find a sequence $\{k_n(1,\cdots,m):n\in{\mathbb N}\}$ such that
$$
\Omega^{\alpha_i}_{r,k_n(1,\cdots,m)}\cap\Omega^{\alpha_j}_{r,k_n(1,\cdots,m)}=\emptyset\s\mbox{for any $i\not= j$.}
$$
Thus,
$$
\varlimsup_{k\rightarrow+\infty}E(u_k,D)\geq \lim_{n\rightarrow+\infty}\sum^m_{i=1} E(u_{k_n(1,\cdots,m)},\Omega_{r,k_n(1,\cdots,m)})\geq m(\tau-\epsilon),
$$
Because $\epsilon>0$ is arbitrary, we are done. \endproof

\begin{defi}\emph{
We say the sequence {\it $\{u_k\}$ has $m$ bubbles }if $\{u_k\}$ has $m$ essentially different blowup sequences and no subsequence of $\{u_k\}$ has more than $m$ essentially different blowup sequences.}
\end{defi}

Although arising from essentially different blowup sequences, some of the $m$ bubbles may be essentially same.

A consequence of Lemma \ref{different} and Lemma \ref{same} is that $\{u_k\}$ has only finitely many bubbles, provided $\sup_kE(u_k,D)<+\infty$. The number of essentially different blowup sequences only depends on the upper bound of the energies.

Now, we assume $\{u_k\}$ has $m$ essentially different blowup sequences $(x_k^1,r_k^1)$, 
$\cdots$, $(x_k^m,r_k^m)$  at $p$, where $m$ depends on $p\in\CS:=\CS(\{u_k\})$. Note that any subsequence $\{u_{k_n}\}$ of $\{u_k\}$ has the same blowup set ${\mathcal S}$ and $\{(x^\alpha_{k_n},r^\alpha_{k_n})\}$ are essentially different  blowup sequences of $\{u_{k_n}\}$. Therefore, to describe the bubble tree convergent it is convenience for us to make the following assumption on the sequence of harmonic maps:

For any $\alpha\not=\beta$, whenever
$\frac{r_k^\alpha}{r_k^\beta}\rightarrow 0$ for any $\alpha\not=\beta$, we always assume, by selecting a subsequence if needed, that \eqref{differentbubble} holds and 
\begin{equation}\label{bubble2}
\mbox{either}\s\frac{|x_k^\alpha-x_k^\beta|}{r_k^\beta}\rightarrow+\infty,\s
\mbox{or}\s\frac{x_k^\alpha-x_k^\beta}{r_k^\beta}\s \mbox{converges.}
\end{equation}
A diagonal sequence procedure ensures the above holds simultaneously for all $1\leq\alpha,\beta\leq m$ and we will still use the index $k$ for the subsequence.

We now build a hierarchy for the bubbles. Put
$$
v_k^\alpha(x)=u_k(x_k^\alpha+r_k^\alpha x)\rightharpoonup v^\alpha,\s\alpha=1,\cdots, m.
$$

\begin{defi}\label{less}{\em
For two essentially different blowup sequences $\{(x^\alpha_k,r^\alpha_k)\}, \{(x^\beta_k,r^\beta_k)\}$, we say $(x_k^\alpha,r_k^\alpha)<(x_k^\beta,r_k^\beta)$, if
$\frac{r_k^\alpha}{r_k^\beta}\rightarrow 0$ and
$\frac{x_k^\alpha-x_k^\beta}{r_k^\beta}$ converges as $k\to\infty$. }
\end{defi}

\begin{lem}
Under the assumption (\ref{bubble2}), it holds
\begin{equation}\label{equiv}
(x_k^\alpha,r_k^\alpha)<(x_k^\beta,r_k^\beta)
\Longleftrightarrow D_{Rr_k^\alpha}(x_k^\alpha)\subset D_{Rr_k^\beta}(x_k^\beta),\s\mbox{for some $R$ and all large $k$.}
\end{equation}
\end{lem}

\proof Assume $(x_k^\alpha,r_k^\alpha)<(x_k^\beta,r_k^\beta)$.  Take 
$$
R=\left|\lim_{k\to\infty}\frac{x^\alpha_k-x^\beta_k}{r^\beta_k}\right| +1.
$$ 
For any $x\in D_{Rr^\alpha_k}(x^\alpha_k)$, 
$$
\frac{|x-x^\beta_k|}{r^\beta_k} \leq \frac{r^\alpha_k}{r^\beta_k}\,\frac{|x-x^\alpha_k|}{r^\alpha_k}+\frac{|x^\alpha_k-x^\beta_k|}{r^\beta_k}
\leq \frac{r^\alpha_k}{r^\beta_k}\ R + R -\frac{1}{2} \leq R\s \mbox{for $k$ large}
$$
thus  $x\in D_{Rr^\beta_k}(x^\beta_k)$. The other direction holds true as $\frac{x^\alpha_k-x^\beta_k}{r^\beta_k}$ is uniformly bounded hence it converges as $k\to\infty$ by (\ref{bubble2}).  \endproof

\begin{lem}\label{still seen}
If $(x_k^\alpha,r_k^\alpha)<(x_k^\beta,r_k^\beta)$,
then the bubble $v^\alpha$ of $\{u_k\}$ at 0 is also a bubble of $\{v_k^\beta\}$
at $x_\beta^\alpha= \lim_{k}\frac{x_k^\alpha-x_k^\beta}{r_k^\beta}$.
\end{lem}

\proof We have
$$
v_k^\alpha(x)=u_k(x_k^\alpha+r_k^\alpha x)=u_k\left(x_k^\beta+r_k^\beta
\left(\frac{r_k^\alpha}{r_k^\beta}x+\frac{x_k^\alpha-x_k^\beta}{r_k^\beta}\right)\right)
=v_k^\beta \left(
\frac{r_k^\alpha}{r_k^\beta}x+\frac{x_k^\alpha-x_k^\beta}{r_k^\beta}\right).
$$
Since $(x_k^\alpha,r_k^\alpha)<(x_k^\beta,r_k^\beta)$,
we have 
$$
\lambda_k:=\frac{r_k^\alpha}{r_k^\beta}\rightarrow 0, \s
x_{\alpha\beta,k}:=\frac{x_k^\alpha-x_k^\beta}{r_k^\beta}\rightarrow x_{\alpha\beta},
$$
 and $v^\beta_k(\lambda_k x +x_{\alpha\beta,k})$ converges as the left hand side $v^\alpha_k$ converges. Hence
$(x_{\alpha\beta,k},\lambda_k)$ is a blowup sequence of
$\{v_k^\beta\}$  at $x_{\alpha\beta}$, and $v^\alpha$ is a bubble of $\{v_k^\beta\}$. \endproof

Intuitively, Lemma \ref{still seen} says that a bubble arises from the ``lower" side in the relation $``<"$ is still captured from the ``upper" side, so it sits on an ``upper" level in the bubble tree. 
If the sequence $v_k^\alpha$ has no concentration points, then
$v^\alpha$ must be at the top of the bubble tree. 
We now make it more precise.

\begin{defi}\label{ontop}{\em
A blowup sequence $(x_k^\alpha,r_k^\alpha)$ is said to be {\it right on top of} another blowup sequence $(x_k^\beta,r_k^\beta)$, if
$(x_k^\alpha,r_k^\alpha)<(x_k^\beta,r_k^\beta)$ and there is no blowup sequence
$(x_k^\gamma,r_k^\gamma)$ that is essentially different from $(x^\alpha_k,r^\alpha_k)$ and $(x^\beta_k,r^\beta_k)$, such that
$(x_k^\alpha,r_k^\alpha)<(x_k^\gamma,r_k^\gamma)<(x_k^\beta,r_k^\beta).$}
\end{defi}

The bubble $v^1$ is rooted at a point on the bubble $v^2$ if $(x^1_k,r^1_k)$ is right on top of $(x^2_k,r^2_k)$. 

\begin{figure}[!ht]
\includegraphics{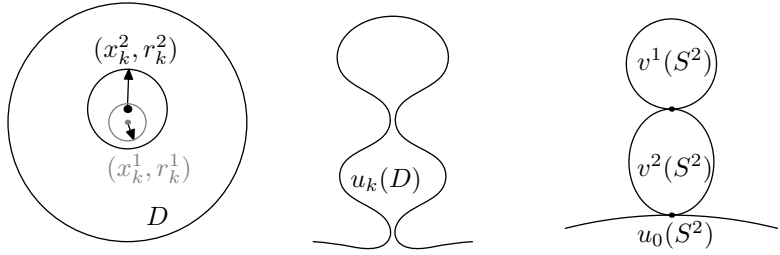}
\caption{Example: $(x_k^1,r_k^1)$ is right on top of $(x_k^2,r_k^2)$.}
\end{figure}

\begin{lem}\label{empty}
If $(x_k^\alpha,r_k^\alpha)$ and $(x_k^\beta,r_k^\beta)$ are essentially different and are both right on top of $(x_k^\gamma,r_k^\gamma)$, then for any fixed $R$, there holds
$$
D_{Rr_k^\alpha}(x_k^\alpha)\cap D_{Rr_k^\beta}(x_k^\beta)=\emptyset
$$
and sufficiently large $k$.
\end{lem}

\proof Assume the result were not true for some $R_0$. Then there would be a subsequence $k_n,n\in{\mathbb N}$ wiere
$$
\frac{|x_{k_n}^\alpha-x_{k_n}^\beta|}{r_{k_n}^\alpha+r_{k_n}^\beta}\leq R_0.
$$
Since $(x_k^\alpha,r_k^\alpha)$ and $(x_k^\beta,r_k^\beta)$ are essentially different
blowup sequence, from (\ref{essen}) we must then have
$\frac{r_{k}^\alpha}{r_{k}^\beta}\rightarrow0$ or
$\frac{r_{k}^\beta}{r_{k}^\alpha}\rightarrow 0$. Without losing generality,
we assume  $\frac{r_{k}^\alpha}{r_{k}^\beta}\rightarrow 0$. By assumption (\ref{bubble2}),  $\frac{x^\alpha_{k}-x^\beta_{k}}{r^\beta_{k}}$ converges, hence $(x_{k}^\alpha,r_{k}^\alpha)<(x_{k}^\beta,r_{k}^\beta)$ by definition. Thus
$(x_k^\alpha,r_k^\alpha)<(x_k^\beta,r_k^\beta)<(x_k^\gamma,r_k^\gamma)$ which contradicts the fact that $(x_k^\alpha,r_k^\alpha)$ is
right on top of $(x_k^\gamma,r_k^\gamma)$.
\endproof

\begin{figure}[!ht]
\includegraphics{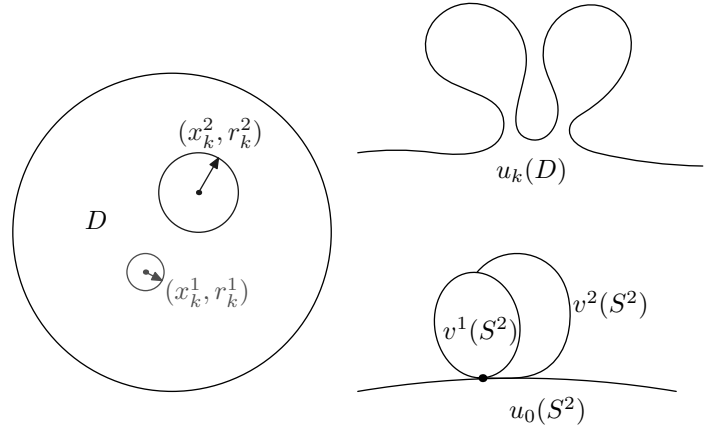}
\caption{Example: $(x_k^1,r_k^1)$ is not right on top of  $(x_k^2,r_k^2)$
and vice versa.}
\end{figure}

\subsection{\bf The neck region} 

To describe the bubble tree convergence of a harmonic map sequence, we need
some notations. For convenience, we call $(x^0_k,r^0_k):=(0,1)$ to be the 0th blow up sequence,
and $u_0$ be the 0th bubble. Then we may say $(x_k^\alpha,r_k^\alpha)<(0,1)$
for any $\alpha\geq 1$.
For a blowup sequence $(x_k^\alpha,r_k^\alpha)$ of $\{u_k\}$ at the center of $D$, let
$\mathcal{S}(\{v_k^\alpha\})$ be the finite set of blowup points of $\{v^\alpha_k\}$, here we recall $v^\alpha_k(x)=u_k(x^\alpha_k+r^\alpha_k x)$. Since there are only finitely many essentially different blowup sequences, there exists $r_0$ such that for $0<r<r_0$ 
\begin{equation}\label{U}
U_{r}^\alpha=D_\frac{1}{r}(0)\backslash 
\bigcup_{p\in \mathcal{S}(\{v_k^\alpha\})} {D_{r}(p)}
\end{equation}
 is  a multiply connected collar and $v_k^\alpha$ converges smoothly on
$U_{r}^\alpha$ when $r$ is fixed. For large $k$, the following region is in $D$: 
\begin{equation}\label{mathcalU}
\Omega_{r,k}^\alpha=r_k^\alpha U_{r}^\alpha+x_k^\alpha.
\end{equation}

For a blowup sequence $(x^\alpha_k,r^\alpha_k)$, let $(x_k^{\sigma_j},r_k^{\sigma_j}), j=1,\cdots,s(\alpha)$, be all of the blowup sequences of $\{u_k\}$ which are right on top of $(x^\alpha_k,r^\alpha_k)$ with
\begin{equation}\label{p}
\lim_{k\to\infty}\frac{x^\alpha_k-x_k^{\sigma_j}}{r_k^\alpha}= p.
\end{equation}
Lemma \ref{still seen} asserts $p\in \CS(\{v_k^\alpha\})$.
Lemma \ref{empty} and (\ref{equiv}) imply
\begin{equation}\label{neck}
N_{r,k}^\alpha(p):= \left(r_k^{\alpha}{D_{r}(p)}+x_k^{\alpha}\right)\backslash
\bigcup_{j=1}^{s(\alpha)}D_{\frac{r_k^{\sigma_j}}{r}}(x_k^{\sigma_j})
\end{equation}
is a multiply connected collar. For each $k$, we call   
$N_{r,k}^\alpha=\bigcup_p N^\alpha_{r,k}(p)$ the {\it neck region for the blowup sequence $(x^\alpha_j,r^\alpha_j)$ at $0\in{\mathcal S}(\{u_j\})$}, where $p$ which satisfies \eqref{p}.  Similarly, we define the neck region for a blowup sequence which is right on top of $(x^\alpha_k,r^\alpha_k)$, and so on. Note that, according to our convention, the neck region for the 0th blowup sequence $(0,1)$ at a blowup point $0\in {\mathcal S}(\{u_k\})$ is 
$D_r(q)\backslash \cup D_{r^\alpha_k/r} (x_k^\alpha)$, where the union is taken over all essentially different blowup sequences $(x^\alpha_k, r^\alpha_k)$ of $\{u_k\}$ at $0$ that are right on top of $(0,1)$.

We will show below that $N_{r,k}$ is a finite union of disjoint multiply connected collars for small $r$ and large $k$. For simplicity of notations, write $\CS^\alpha$ for $\CS(\{v^\alpha_k\})$.

We choose $r_0$ to satisfy: 
$$r_0< \frac{1}{2} \min_{p\not=q}|p-q|,\s \mbox{ for $p,q\in\cup_\alpha \CS^\alpha$} 
$$ 
and 
$$
r^{-1}_0-r_0> 1+\max_{p\in\cup_\alpha\CS^\alpha}|p|.
$$ 
Then for any fixed $0<r<r_0$, it holds that for any  $\beta$ and $q\in \cup_\alpha\CS^\alpha$,
\begin{equation}\label{1}
N_{r,k}^{\beta}(q)
\subset D_{\frac{r_k^{\beta}}{r}}(x_k^{\beta})
\end{equation}
since $r^{-1}-r>|q|$ implies  
$$
x_k^{\beta}+r_k^{\beta} D_r(q)\subset D_{r_k^{\beta}(|p|+r)}(x_k^{\beta}).
$$

\begin{pro}\label{nointersection}
 For sufficiently large $k$ and $r<r_0$, it holds 
\begin{enumerate}
\item[(a)] $N^\alpha_{r,k}(p)\cap N^\alpha_{r,k}(q)=\emptyset$ for any $p\not=q$.
\item[(b)] $N^\alpha_{r,k}(p)\cap N^\beta_{r,k}(q)=\emptyset$ for any $\alpha\not=\beta$.
\end{enumerate}
 \end{pro}

\proof 

(a) By definition, $N^\alpha_{r,k}(p)\subset x^\alpha_k+r^\alpha_k D_r(p), N^\alpha_{r,k}(q)\subset x^\alpha_k+r^\alpha_k D_r(q)$, and $r<\frac{1}{2}|p-q|$. Thus $N^\alpha_{r,k}(p)\cap N^\alpha_{r,k}(q)=\emptyset$. 

(b) For any $\alpha\not=\beta$, since $(x^\alpha_k,r^\alpha_k), (x^\beta_k,r^\beta_k)$ are essentially different and we have assumed (\ref{bubble2}), it holds either
\begin{enumerate}
\item  ${\frac{x^\alpha_k-x^\beta_k}{r^\alpha_k+r^\beta_k}\to\infty}$ and $c_1< \frac{r^\alpha_k}{r^\beta_k}<c_2$ for some positive constant $c_1,c_2$; or
\item ${\frac{r^\alpha_k}{r^\beta_k}\to 0}$ and ${\frac{x^\alpha_k-x_k^\beta}{r^\alpha_k+r^\beta_k}}$ converges.
\end{enumerate}

If $\left(r_k^{\alpha}{D_{r}(p)}+x_k^{\alpha}\right) \cap (r_k^{\beta}{D_{r}(q)}+x_k^{\beta})=\emptyset$ for all large $k$, then (b) holds since $N^\alpha_{r,k}(p)\subset r^\alpha_k D_r(p)+x^\alpha_k$ and $ N^\beta_{r,k}(q)\subset r^\beta_k D_r(q)+x^\beta_k$. 

If there are $x_k\in \left(r_k^{\alpha}{D_{r}(p)}+x_k^{\alpha}\right) \cap (r_k^{\beta}{D_{r}(q)}+x_k^{\beta})$ for infinitely many $k$, then
$\frac{x_k-x^\alpha_k}{r^\alpha_k} \in {D_{r}(p)}$ and $ \frac{x_k-x^\beta_k}{r^\beta_k}\in {D_{r}(q)}$, here we still use index $k$ for simplicity. Therefore (1) cannot happen as
$$
|x^\alpha_k-x^\beta_k| \leq |x_k-x^\alpha_{k}|+|x_k-x^\beta_k| \leq (r+|p|)r^\alpha_k+(r+|q|)r^\beta_k\leq (r+|p|+|q|)(r^\alpha_k+r^\beta_k).
$$ 
If (2) happens, then $(x^\alpha_k,r^\alpha_k)<(x^\beta_k,r^\beta_k)$ and
 $\frac{x_k^\alpha-x_k^\beta}{r_k^\beta}$ converges to some blowup point $\bar{q}$ of the sequence $u_k(x^\beta_k+r^\beta_kx)$.

Case 1.   $\bar{q}\neq q$. Since $r<\frac{|\bar{q}-q|}{2}$, it follows 
\begin{equation}\label{2}
(x_k^\beta+r_k^\beta D_r(\bar{q}))\cap
 (x_k^\beta+r_k^\beta D_r(q))=\emptyset.
\end{equation}
Given $x\in D_\frac{r_k^\alpha}{r}(x_k^\alpha)$, we can write
$x=x_k^\alpha+r_k^\alpha y$ for some $y\in D_\frac{1}{r}(0)$. Then
$$
x=x_k^\beta+r_k^\beta \left(\frac{x_k^\alpha-x_k^\beta}{r_k^\beta}+\frac{r_k^\alpha}{r_k^\beta}y\right).
$$
Recalling that $\frac{x_k^\alpha-x_k^\beta}{r_k^\beta}\rightarrow \bar{q}$ and
$\frac{r_k^\alpha}{r_k^\beta}\rightarrow 0$, we have 
$\frac{x_k^\alpha-x_k^\beta}{r_k^\beta}
+\frac{r_k^\alpha}{r_k^\beta}y\in D_r(\bar{q})$ for a fixed $r$ and sufficiently large $k$.
Thus,
$$
D_\frac{r_k^\alpha}{r}(x_k^\alpha)\subset x_k^\beta+r_k^\beta D_{r}(\bar{q}), 
$$
which can be combined with  \eqref{1} to yield
$$
N^\alpha_{r,k}(p)\subset D_{\frac{r^\alpha_k}{r}}(x^\alpha_k)\subset x^\beta_k+r^\beta_kD_r(\bar{q}). 
$$
By definition, 
$$
N^\beta_{r,k}(q)\subset x^\beta_k+r^\beta_kD_r(q).
$$ 
Now it follows from  \eqref{2} that 
$$
N_{r,k}^\alpha(p)\cap N_{r,k}^\beta(q)=\emptyset.
$$

Case 2. $\bar{q}=q$ and $(x_k^\alpha,r_k^\alpha)$ is right on top of  $(x_k^\beta,r_k^\beta)$. In this case, the definition of $N^\beta_{r,k}(q)$ directly leads to $N^\beta_{r,k}(q)\cap 
D_{{r^\alpha_k}/{r}}(x^\alpha_k)=\emptyset$.  Note  
$N^\alpha_{r,k}(p) \subset D_{{r^\alpha_k}/{r}}(x^\alpha_k)$ by \eqref{1}. It then follows 
$N_{r,k}^\alpha(p)\cap N_{r,k}^\beta(p)=\emptyset$.

Case 3. $\bar{q}=q$ and there is $(x_k^\gamma,r_k^\gamma)$ with $(x_k^\alpha,r_k^\alpha)<(x_k^\gamma,r_k^\gamma)<(x^\beta_k,r^\beta_k)$. As there are only finitely many blowup sequences, we may assume that $(x^\gamma_k,r^\gamma_k)$ is right on top of $(x^\beta_k,r^\beta_k)$. 
Then we have
$$
\frac{r_k^\alpha}{r_k^\gamma}
\rightarrow 0,\s\frac{r_k^\gamma}{r_k^\beta}\rightarrow 0,\s 
\frac{x_k^\alpha-x_k^\gamma}{r_k^\gamma} \rightarrow q',\s\mbox{for some $q'$}.
$$
Then we have
$$
\frac{x_k^\gamma-x_k^\beta}{r_k^\beta}=\frac{x_k^\alpha-x_k^\beta}{r_k^\beta}
+\frac{x_k^\gamma-x_k^\alpha}{r_k^\gamma}\cdot \frac{r_k^\gamma}{r_k^\beta}\rightarrow \bar{q}.
$$
For any $x\in D_{\frac{r_k^\alpha}{r}}(x_k^\alpha)$, we have for $k$ large 
$$
|x-x_k^\gamma|\leq |x-x_k^\alpha|+|x_k^\alpha-x_k^\gamma|\leq 
\left(\frac{r_k^\alpha}{rr_k^\gamma}+\frac{|x_k^\alpha-x_k^\gamma|}{r_k^\gamma}\right)r_k^\gamma
\leq (1+|q'|)r_k^\gamma.
$$
Since $\frac{1}{r}>|q'|+1$, then we get 
$$D_{\frac{r_k^\alpha}{r}}(x_k^\alpha)\subset
 D_{\frac{r_k^\gamma}{r}}(x_k^\gamma).$$ 
The definition of $N_{r,k}^\beta(q)$ implies 
$N_{r,k}^\beta(q)\cap D_{{r_k^\gamma}/{r}}(x_k^\gamma)
=\emptyset$. By \eqref{1}, $N^\alpha_{r,k}(p)\subset D_{{r^\alpha_k}/{r}}(x^\alpha_k)$, so we conclude  $N_{r,k}^\alpha(p)\cap N_{r,k}^\beta(p)=\emptyset$.

This completes the proof of (b).  \endproof


\subsection{\bf Convergence}

By results in \cite{Parker,Ding-Tian,Qing,Qing-Tian}, a sequence of harmonic maps $\{u_k\}$ with uniformly bounded energy, or more generally, a Palais-Smale sequences for the energy functional, with uniformly $L^2$-bounded tension fields,  converges in the sense of the bubble tree as follows
$$
u_k\rightarrow u_\infty:= (u_0;\s v^1,\s \cdots,\s v^m),
$$
where $u_0$ is a harmonic map from $D$ and each $v^i$ is a bubble \cite{S-U}, with the following properties: 

\begin{enumerate}

\item[C1.] $u_k\rightharpoonup u_0$ in $W^{1,2}(D)$.

\item[C2.] There exist essentially different blowup sequences, i.e. satisfy
\eqref{bubble1},
$(x_k^\alpha,r_k^\alpha)$, $\alpha =1,2,\cdots,m$, such that
$$u_k^\alpha(x)=u_k(x_k^\alpha+r_k^\alpha x)\rightharpoonup v^\alpha,\s
\mbox{in}\s W^{1,2}_{loc}(\R^2).$$

\item[C3.] When $\frac{r_k^\alpha}{r_k^\beta}\rightarrow 0$, \eqref{bubble2}
holds.

\item[C4.] For any fixed $r$, $u_k^\alpha$ converges smoothly to $v^\alpha$
on ${U}_{r}^\alpha$, where $U_{r}^\alpha$ is defined by \eqref{U}.

\item[C5.] Energy identity: $\lim\limits_{k\rightarrow+\infty}E(u_k,D_\frac{1}{2})=
E(u_0,D_\frac{1}{2})+\sum\limits E(v^\alpha)$. 

\item[C6.] No neck: for any $\alpha$, we have
$$\lim_{r\rightarrow 0}\lim_{k
\rightarrow+\infty}\sup_{x\in N_{k,r}^\alpha(p)}
|u_k(x)-v^\alpha(p)|=0
.$$\\
\end{enumerate}
The bubble tree limit $u_\infty$ can be viewed as a continuous map from a bubble tree formed by $D$ and $m$ copies of $S^2$.


\section{Proof of Theorem \ref{str}} 

In this section, we first construct normalizations for harmonic maps from a single $S^2$, and use them to show that there are only finitely many homotopy classes below a given energy level. This of course follows from Parker's bubble tree convergence theorem \cite{Parker}, without using the normalizations. 
Next, we construct normalizations of stratified spheres from the ones for a single $S^2$ and prove Theorem \ref{str}.

\subsection{Homotopy class of  harmonic spheres with bounded energy}\label{B.T.}

Suppose that $\{u_k\}$ is a sequence of harmonic maps from $S^2$ into $N$ with uniform bounded
energies which
converges to $(u_0$; $v^1$, $\cdots$, $v^m)$ in the sense of bubble tree, where $v^1,\cdots,v^m$ are all non-trivial harmonic maps from $S^2$ to $N$.
We will construct a stratified sphere $\S_\infty$ and
put  $u_0$, $v^1$, $\cdots$, $v^m$ together to obtain
a harmonic map from $\S_\infty$ to $N$. We follow the following steps to construct $\S_\infty$:

\vspace{.25cm}

Step 1.  Around each blowup point of $\{u_k\}$, call it $p_0$, take an isothermal coordinate system $(D,x)$ around $p_0$ such that $p_0=0$.
Let $v^1$, $v^2$, $\cdots$, $v^{l_1}$ be all of the
bubbles which are right on top of $u_0$, i.e. the corresponding blowup sequences are right on top of that of $u_0$, which, in this initial step, is $(0,1)$ in our convention. Each $v^\alpha$ is a harmonic map from a standard sphere $S_i$ to $N$. We identify the south pole of $S_i$ with $p_0\in S_0=S^2$.

\vspace{.25cm}

Step 2. Repeat Step 1 at

$(\mbox{a}_1)$ each concentration point of $\{v_k^j \}$
where  $v_k^j(y)=u_k(y_k^j+\rho_k^j y)$,

$(\mbox{a}_2)$ each concentration point of $\{w^s_k\}$ where $w^s(z)=v^s_k(z^s_k+\sigma^s_k z)$

$(\mbox{a}_3)$ similar to the above two cases, at all higher levels in the bubble tree.

\vspace{.25cm}

In the end, we get a stratified sphere $\S_\infty$ and a harmonic map
$u_\infty:\S_\infty\rightarrow N$.

\vspace{.25cm}

Next we construct a normalization $\phi_{r,k}$ of $\S_\infty$ such that
$u_\infty\circ\phi_{r,k}$ is very close to $u_k$.


Let $q_1, \cdots, q_l$ be the blowup points of $\{u_k\}$. 
Now, on an isothermal coordinate system around each concentration point $q_i=0$ (we may assume that $q_i$ is the only blowup point of $\{u_k\}$ in the local chart),
we define a continuous map $\phi_{r,k}:S^2\to\S_\infty$   as follows. 

At first level: outside the disks $D_r(q_i)$'s:
\begin{enumerate}

\item[1)] $\phi_{r,k}$ is the identity map on $S^2\backslash \cup_i D_{2r}(q_i)$;

\item[2)] $\phi_{r,k}$ is bijective from $D_{2r}(q_i)\backslash D_r(q_i)$ onto
$D_{2r}(q_i)\backslash\{q_i\}$, for each $i=1,\cdots,l$;

\item[3)] $\phi_{r,k}(\partial D_r(q_i))=q_i$;

\end{enumerate}

At the second level:  inside $D_r(q_i)$'s:

We first recall from Proposition \ref{nointersection} for all large $k$ and $r<r_0$ that $N^\alpha_{r,k}(p)$'s are disjoint and lie inside $D_r(q_i)$, and 
$N^\alpha_{r,k}(p) =\left( r^\alpha_k D_r(p)+x^\alpha_k \right)\backslash \cup D_{{r^{\sigma_j}}/{r}}(x^{\sigma_j}_k)$. Select those $\alpha$ and $p$ such that 
the corresponding disks $D_k^i(p_i)= r^\alpha_k D_r(p)+x^\alpha_k$ are disjoint. Define $\phi_{r,k}$ bijectively from $D_r(q_i)\backslash \cup_i D_k^i(p_i)$
onto $D_r(q_i)\backslash\{ p_1,...,p_s\}$ and maps $\partial D^i_k(p_i)$ to $p_i$. Further, define $\phi_{r,k}$ on $N^{\alpha_i}_{r,k}(p_i)$ by sending every point to $v^{\alpha_i}(p_i)$. 

Third and higher levels: on the disks $D_{{r^{\sigma_j}}/{r}}(x^{\sigma_j}_k)$ which are inside $\cup_i D^i_k(p_i)$. Repeat the second step. 

Finally, on the remaining disks which correspond to the bubbles at the highest level, i.e. there are no blowup sequences right on top of those in the previous steps, note the boundary of each of these disks is mapped to a point in the previous step, now define $\phi_{r,k}$ be a diffeomorphism from the disk to $S^2$ with the point deleted. This completes the definition of $\phi_{r,k}$. 


As there is no neck, we have
\begin{equation}\label{noneck}
\lim_{r\rightarrow 0}\lim_{k\rightarrow+\infty}\sup_{x\in S^2}
d_N(u_\infty\circ\phi_{r,k}(x),u_k(x))=0.
\end{equation}
Thus, $u_\infty\circ\phi_{r,k}\sim u_k$. It is easy to check that $\phi_{r,k}$ is a
normalization of $\S_\infty$ . Then by definition $u_\infty\sim u_k$. Hence,
\begin{lem}
Let $u_k$ be a sequence of harmonic maps from
the sphere $S^2$ to $N$ with uniformly bounded energies. Then after passing
to a subsequence, all $u_{k}$ are  in the same homotopy class.
\end{lem}

\subsection{Bubble tree convergence of harmonic maps from stratified spheres}

Let $\{x_k\}$ be a sequence in $S^2$. We say $x_0\in \S_\infty$ is the {\it $\phi$-limit of $\{x_k\}$}, if
$$x_0=\lim_{r\rightarrow0}\lim_{k\rightarrow+\infty}\phi_{r,k}(x_k).$$
By the no-neck result, we have
\begin{equation}\label{prelimit}
\lim_{k\rightarrow+\infty}u_k(x_k)=u_\infty(x_0).
\end{equation}
It is easy to check that for any sequence $x_k\in S^2$, we
can find a subsequence which has a $\phi$-limit .\\

The bubble tree convergence of $\{u_k\}$ can also be described via $\phi_{r,k}$.
Let $P$ be the set of singular  points of $\S_\infty$. Let
$S_\delta=\S_\infty\backslash\cup_{p\in P}B_\delta(p)$.
Clearly, for any $\delta>0$, $F_{r,k}=\phi_{r,k}^{-1}$
is well defined on $S_\delta$ when $r$ is sufficiently small. Moreover  $u_k\circ F_{r,k}$
is harmonic and converges smoothly. If $V_{\delta,k}$
is a connected component of $S^2\backslash F_{r,k}(S_\delta)$, then for any $p\in P$
$$\lim_{\delta\rightarrow 0}\lim_{k\rightarrow+\infty}\underset{V_{\delta,k}}{\mbox{osc}} \,u_k=0.$$


Let $\S_k=S_1\cup S_2\cup\cdots \cup S_{m(k)}$ and $P_k$ be the set of singular points of $\S_k$, where $S_i$'s are smooth 2-spheres. Note that $P_i\neq P_j$ may happen.

Let $u_k$ be a harmonic map from $\S_k$ into $N$, i.e. $u_k$ is continuous on $\S_k$
and harmonic on each component of $\S_k$. We may assume $u_k|_{S_i}$
converges to $(u_\infty^i,v^{i,1},\cdots,v^{i,l_i})$ in the sense of bubble tree, where
$v^{i,j}$ is not trivial. Then
we get a stratified sphere $\S^i_\infty$, and a harmonic map
$u_\infty^i:\S^i_\infty\rightarrow N$.

Now, let $p_k^1\in P_k$ be the singular point which connects $S_1$ and $S_2$. By passing to subsequences if needed, we may assume $p^1$ is the $\phi$-limit  of $p_k^1$ on $\S^1_\infty$ and $p^2$ is the $\phi$-limit  of $p_k^1$ on $\S^2_\infty$. By the no-neck result,  $u_\infty^1(p^1)=u_\infty^2(p^2)$.
We then identify $p^1$ and $p^2$ to get a stratified sphere $\S^{12}_\infty$, and  a harmonic map
from $\S^{12}_\infty$ into $N$.
Repeating this construction for $\S^i_\infty,\S^j_\infty$ whenever $S_i,S_j$  are joint at a singular point, we get a stratified sphere $\S_\infty$ and a harmonic map
$u_\infty:\S_\infty\to N$. We may view $\S_\infty$ as the union of $\S_\infty^i$.

\subsection{Proof of Theorem \ref{str}:} By the gap theorem for 
harmonic maps from $S^2$ \cite{S-U}, the number of components of $\S_k$ is uniformly bounded in $k$. Therefore there are only finitely many homeomorphic types. It suffices for us to assume $\S_k$ are of the same homeomorphic type.

Let  $u_k:\S_k\to N$ be harmonic maps with the uniformly bounded
energies. 
We assume $u_k$ converges in the bubble tree sense to a harmonic map
$u_\infty:\S_\infty\rightarrow N$ as we described in the previous subsection.
Moreover, we may set $\S_k=S_1\cup\cdots \cup S_{m}$,
$\S_\infty=\S_{\infty}^1\cup \cdots \cup \S_\infty^m$,
$u_k^i=u_k|_{S_i}$, such that $u_k^i$ converges in the sense
of bubble tree to $u_\infty^i:\S_\infty^i\rightarrow N$, where
$u_\infty^i=u_\infty|_{\S_\infty^i}$.

Let $\phi_{r,k}^i$ be the normalization of $\S_\infty^i$, which is defined in subsection \ref{B.T.}. To use Lemma \ref{homotopy},
we need to combine $\phi_{r,k}^i$ over $i$ to produce a continuous map from $\S_k$ to $\S_\infty$.
To achieve this,  we need to modify $\phi_{r,k}^i$ on neighborhoods of the singular points.

Let $x_k\in S_i\cap P_k$, such that $\{x_k\}$ has a $\phi$-limit $x_0\in\S_\infty^i$.
We have two cases:

Case 1. $x_0$
is a singular point of $\S_\infty^i$. In this case, we can suitably choose $r$, so that
$\phi_{r,k}^i(x_k)=x_0$. 
Set $\varphi_k^i = \mbox{the identity map}$ of $\S^i_\infty$.

Case 2. $x_0$ is not a singular point of $\S_\infty^i$. Choose a disk $B_\epsilon(x_0)\subset
\S_\infty^i$, such that $\phi_{r,k}^i(x_k)\in B_\epsilon(x_0)$ for large $k$ and $B_\epsilon(x_0)$ contains no singular points of $\S^i_\infty$. Let $\varphi_k^i$ be a homeomorphism
 of $\S_\infty^i$ onto itself such that $\varphi_k^i\circ\phi_{r,k}^i(x_k)=x_0$ and
$\varphi_k|_{\S_\infty^i\backslash B_\epsilon(x_0)}$ is the identity map.
Thus, $\varphi_k^i\circ\phi_{r,k}^i$ is a normalization of $\S_\infty^i$ with
$\varphi_k^i \circ\phi_{r,k}^i(x_k)=x_0$.

Thus, $f_k:\S_k\to \S_\infty$, defined by $f_k=\varphi_k^i\circ\phi_{r,k}^i$ on each $S_i$,
is a continuous map and $f_k|_{S_i}$ is a normalization of $\S_\infty^i$.
By Lemma \ref{homotopy}, $u_\infty\sim u_\infty\circ f_k$.

Furthermore, in Case 2, we may suitably choose
$\varphi_k^i$, as $x_0$ is a $\phi$-limit of $\{x_k\}$, such that $\varphi_k^i$ converges to the identity map as $k\rightarrow+\infty$ (this is always true for Case 1). Note that
$$
\sup_{x\in (\phi_{r,k}^i)^{-1}(B_\epsilon(x_0))}
d_N(u_\infty\circ\varphi_k^i\circ\phi_{r,k}^i(x),u_k(x))=\sup_{y\in B_\epsilon(x_0)}
d_N\left(u_\infty\circ\varphi_k^i(y),u_k\circ(\phi_{r,k}^i)^{-1}(y)\right),
$$
and recall that $u_k\circ(\phi_{r,k}^i)^{-1}$ converges smoothly on $B_\epsilon(x_0)$. Then for any
$\sigma$ we can suitably choose $r$, such that
\begin{equation}\label{epsilon'}
\sup_{x\in (\phi_{r,k}^i)^{-1}(B_\epsilon(x_0))}
d_N \left(u_\infty \circ \varphi_k^i\circ\phi_{r,k}^i(x),u_k(x) \right)<\sigma.
\end{equation}
By \eqref{noneck},  we may assume
$$
d_N \left( u_\infty^i\circ \phi_{r,k}^i(x),u_k(x) \right)<\sigma.
$$
Then
$$
d_N \left( u_\infty \circ f_k(x),u_k(x) \right)<\sigma.
$$
Thus, $u_\infty \circ f_k \sim u_k$ and it follows $u_k\sim u_\infty$ since we have shown $u_\infty\sim u_\infty\circ f_k$.
\endproof

\section{Homotopy classes of $S^2$ in 3-manifolds}

In this section, we will construct a family of smooth maps from $S^2$ into some 3-dimensional manifold which are not homotopic to each other while their energies  are uniformly bounded above. The elements $\gamma_i$ in the fundamental group $\pi_1(M)$ of $M$ which yield distinct homotopy classes are produced in Proposition \ref{h.sequence} in Appendix.

\begin{thm}\label{3-manifold}
Let  $M=M_1\# M_2$, where $M_1\neq S^3$ and $M_2\neq S^3$ or $\R\P^3$.  
Then there is a sequence of smooth mappings 
$u_k:S^2\rightarrow M$, such that $\sup_kE(u_k)<+\infty$ and
$u_i$ is not homotopic to $u_j$ for any $i\neq j$.
\end{thm}

\proof Let $\pi:\widetilde{M}\to M$ be the universal cover of $M$ and let $\tilde{u}$, $\gamma_i$
be the ones in Proposition \ref{h.sequence}.

Let $P$ and $N$ be the south pole and the north pole of
$S^2$ respectively.
Let $\tilde{u}'$ be a map from $S^2$ to $\widetilde{M}$
which is homotopic  to $\tilde{u}$, such that $\tilde{u}'$ is a constant in a neighborhood
of $P$ and a neighborhood of $N$ respectively.
Let $\Phi=\tilde{u}'\circ \Pi^{-1}:\C\rightarrow \widetilde{M}$, where $\Pi$
is the stereographic projection from $P$.

We assume $\gamma_0=1$. Let $\beta_i:[0,1]\rightarrow \hat{M}$ be an embedded curve with $\beta_i(0)=
\tilde{u}'(P)$ and $\beta_i(1)=\gamma_i\circ \tilde{u}'(N)$. We define $v_k:\C\rightarrow M$
as follows:
$$v_k=\left\{\begin{array}{ll}
\Phi(x), &|x|\geq\delta,\\
\beta_k\left(\frac{\log r-\log R\epsilon}{\log\delta-\log
R\epsilon}\right), &R\epsilon<|x|<\delta,\\
\gamma_k\left(\Phi\left(\frac{x}{\epsilon}\right)\right),&\frac{|x|}{\epsilon}
\leq R,
\end{array}\right.$$
where $\delta$ be sufficiently small and $R$ be sufficiently large, such that
$\Phi$ is constants on $D_\delta$ and $\C\backslash D_R$.
We have
$$\begin{array}{lll}
  \dint_{D_\delta\backslash D_{R\epsilon}}|\nabla v_k|^2&=&
     2\pi\dint_{R\epsilon}^\delta\left|\frac{\partial\beta_k}{\partial r}\right|^2rdr\\
     &<&\frac{c\|\dot{\beta_k}\|_{L^\infty}^2}
     {(-\log R\epsilon+\log\delta)^2}
\dint_{R\epsilon}^\delta\frac{dr}{r}
   =\frac{c\|\dot{\beta_k}\|_{L^\infty}^2}{\log\delta-\log{R\epsilon}},
\end{array}$$
$$\dint_{\C\backslash D_\delta}|\nabla v_k|^2 \leq E(\Phi),$$ and
$$\dint_{D_{R\epsilon}}|\nabla v_k|^2 \leq E(\Phi).$$ So, we can find
suitable $\epsilon$, such that
$$E(v_k)\leq 2E(\Phi)+1.$$

Since $v_k\equiv \tilde{u}_k'$ on $\C\backslash D_R$, we can view $v_k$
as a map from $S^2$ into $\widetilde{M}$. We need to prove that $v_i$ is not homotopic
to $v_j$ for any $i\neq j$. By the Hurewicz Theorem, we only need to check that
$[v_i]_{H_2}\neq [v_j]_{H_2}$. Obviously,
$$[v_i]_{H_2}=[\tilde{u}']_{H_2}+[\gamma_i\circ\tilde{u}]_{H_2},\s [v_j]_{H_2}
=[\tilde{u}']_{H_2}+[\sigma_j\circ\tilde{u}]_{H_2}.$$
Since $[\gamma_i\circ\tilde{u}]_{H_2}\neq [\gamma_j\circ\tilde{u}]_{H_2}$, we have
$[v_i]_{H_2}\neq [v_j]_{H_2}$.

Now, we let $u_k=\pi(v_k)$. Then $u_i$ is not homotopic to $u_j$ for any $i\neq j$. Moreover,
it is easy to check $E(u_k)=E(v_k)$.
\endproof

\section{harmonic map flow with finite time blowup}

Theorem \ref{str} can be used to produce harmonic map flows of $S^2$ in closed 3-manifolds which develop singularities in a finite time.

\begin{thm}\label{harmonicblowup}
Let   $M_1\neq S^3$ and $M_2\neq S^3$ or $\R\P^3$ be closed 3-dimensional Riemannian manifolds and $M=M_1\# M_2$. 
Then there exist infinitely many smooth maps $u_k:S^2\to M$ such that $u_i,u_j$ are not homotopic for any $i\neq j$, $\sup_m E(u_k)<\infty$ and 
the harmonic map flow that begins at $u_k$ develops a finite time singularity. 
\end{thm}

\proof
Theorem \ref{3-manifold} guarantees existence of a sequence of mutually nonhomotopic smooth maps  $v_k:S^2\to M$ with uniformly bounded energies, and only finitely many (may be none) of them contain harmonic maps from some stratified spheres to $M$ in their homotopy classes by Theorem \ref{str}. 
Select the $v_k$'s which do not have harmonic maps from stratified sphere to $M$ in their homotopy classes and denote them by $u_k$. If 
the harmonic map flow
$$
\frac{\partial u}{\partial t}=-\tau(u),\s \s u(x, 0)=u_k(x)
$$
admits a longtime solution, then by results
in \cite{Qing-Tian}, as $t\rightarrow+\infty$, $u(x,t)$ will
converge to a harmonic map from a stratified
sphere to $M$ which is homotopic to $u_k$. This contradicts the choice of $u_k$.
\endproof



\section{Mean Curvature Flow with finite time blowup}
It is well-known that a minimal sphere is also a harmonic sphere. In this section, we combine Theorem \ref{str}, Theorem \ref{3-manifold} and the compactness theorem in \cite{C-L} to produce mean curvature flows that blowup in a finite time in closed 3-manifolds.

\begin{thm}\label{MCF}
Let   $M_1\neq S^3$ and $M_2\neq S^3$ or $\R\P^3$ be closed 3-dimensional Riemannian manifolds and $M=M_1\# M_2$. 
Then there exist infinitely many embeddings $w_k:S^2\to M$ such that $w_i,w_j$ are nonhomotopic for any $i\neq j$ and the areas $\mu(w_k(S^2))$ are uniformly bounded, and the mean curvature flow of $S^2$ initiated at $w_k$ develops a singularity in a finite time.
\end{thm}

\proof First, select $u_k$ as in the proof of Theorem \ref{harmonicblowup}. 
Next,  we replace the line $\beta_k$ in $u_k$ by a thin cylinder to get an embedding $v_k$ from $S^2$ with uniformly bounded area and $v_k\sim u_k$. This can be done as follows. Recall, in the proof of Theorem \ref{3-manifold}, $u_k$ maps the annulus $R\epsilon<|x| < \delta$ to the curve $\beta_k$ and $\Phi$ maps the domain $|x| > R_0>\delta$ to a point and the disk $|x|< R^{-1}_0<R\epsilon$ to a point for some large $R_0$ (these two domains correspond to the two small caps in $S^2$ that are mapped to points by $\tilde{u}'$), and  $\Phi$ is diffeomorphic on the complement of these two domains. Extend $\beta_k$ a little at its end points in $M$, and take the thin embedded cylinder in the tubular neighborhood of the extended $\beta_k$ such that it is the circle of radius $r$ in the normal plane at each point of the extended $\beta_k(t)$. We can assume that one boundary circle of the cylinder is $\Phi|_{\{ |x| = R\epsilon+1/k\}}$ pointwise, and the other boundary circle of the cylinder is $\Phi|_{\{|x|= \delta+1/k\}}$ twisted by a local diffeomporphism.  Define $v_k$ by sending each ray $(r,\theta)$ in the annuls $ R\epsilon +1/k < |x| < \delta - 1/k$ to the corresponding curve in the twisted cylinder in $M$, and outside the annulus let $v_k$ be $u_k$.  We can assume $v_k$ is smooth by slight modification in a neighborhood of the boundary circles.

Note that $v_k$ is conformal, in fact isometric, with respect to the pull back metric $v^*h$ on $S^2$ from $M$.  Since there is only one conformal structure on $S^2$, there exists a diffeomorphism $\phi_k:S^2\to S^2$ such that $v^*_kh= \lambda^2 \phi_k^* g_0$, where $g_0$ is the round metric on $S^2$. 
Due to conformality, area of $v_k(S^2)) = E(v_k, v_k^*h)$, and conformal invariance of energy implies 
$ E(v_k, v^*_kh) = E( v_k\circ \phi_k, g_0)$. Set $w_k = v_k\circ \phi_k$. As any diffeomorphim of $S^2$ to $S^2$ is homotopic to the identity map of $S^2$, we have $w_k \sim v_k$. In conclusion, $w_k$'s are conformal embeddings from $S^2$ to $M$ with uniformly bounded energies in the round metric, and $w_i$ is not homotopic to $w_j$ for $i\not= j$.

Assume the mean curvature flow $F:S^2\times [0,T)\to M$ begins at $w_k$ exists for longtime. Write $\Sigma_t$ for the image of $S^2$ at time $t$. Since
$$
\int^{t}_{0}\int_{\Sigma_t}|H(\Sigma_t)|^2d\mu_t = \mu(\Sigma_0)-\mu(\Sigma_t) <\infty
$$
for all $t\in [0,\infty)$, there exists a sequence $t_n\to\infty$ such that
$$
\lim_{n\to\infty}\int_{\Sigma_{t_n}}|H(\Sigma_{t_n})|^2d\mu_{t_n} = 0.
$$
Each surface $\Sigma_{t_n}$ is immersed, and $F(\cdot, t_n): (S^2, F^*(\cdot,t_n)h) \to (M,h)$ is conformal with respect to the pull-back metric from $M$. 
As there is only one conformal class on $S^2$, 
$F(\cdot,t_n)\circ \phi_n:(S^2,g_0)\to (M,h)$ is conformal for some diffeomorphism $\phi_n:S^2\to S^2$. Along the mean curvature flow, the area of the conformal maps $F(\cdot,t_n)\circ\phi_n$ is uniformly bounded

$$
\left| F(\cdot,t_n)\circ\phi_n(S^2) \right| = \left|F(S^2,t_n)\right|\leq \left|F(S^2,0)\right|=E(w_k)<C, 
$$ 
and  so is the Willmore energy
$$
\int_{F(\cdot,t_n)\circ\phi_n(S^2)} |H_{F(\cdot,t_n)\circ\phi_n}|^2 d\,\mu_{F(\cdot,t_n)\circ\phi_n} = \int_{\S_{t_n}} \left| H(\S_{t_n})\right|^2 d\,\mu_{t_n)} \to 0. 
$$

We now recall a special case (namely, domain is $S^2$) of the compactness theorem  in \cite{C-L}: Suppose that $\{f_k\}$ is a sequence of $W^{2,2}$  branched conformal
immersions of $(S^2, g_0)$ in a compact manifold $M$. If
$$
\sup_k\left\{\mu(f_k)+W(f_k)\right\}<+\infty
$$
then  either $\{f_k\}$ converges to a point, or
there is a stratified sphere $\Sigma_\infty$ and a $W^{2,2}$ branched conformal immersion 
$f_\infty:\Sigma_\infty\to M$,
such that a subsequence of
$\{f_k(\Sigma)\}$ converges to  $f_\infty(\Sigma_\infty)$
in the Hausdorff topology, and the area and the Willmore energy satisfy
$$
\mu(f_\infty)=\lim_{k\rightarrow+\infty}\mu(f_k)\s \mbox{and} \s W(f_0)\leq\lim_{k\rightarrow+\infty} W(f_k).
$$

By this compactness theorem, there is a stratified sphere $\Sigma_\infty=\Sigma^1_\infty \cup \cdots \cup \Sigma^k_\infty$, where each $\Sigma^i_\infty$ is a 2-sphere $S^2$, and a $W^{2,2}$ branched conformal immersion $F_\infty:\Sigma_\infty\to M$ with $\Sigma_{t_n}$ converging to $F_\infty(\Sigma_\infty)$ in the sense of bubble tree, therefore $F_\infty$ is homotopic to $w_k$. Moreover, 
$$
\int_{F_\infty(\Sigma_\infty)}|H|^2d\mu \leq \lim_{n\to\infty} \int_{\Sigma_{t_n}}|H|^2d\mu_{t_n} =0.
$$
We see that $F_\infty|_{\Sigma_{\infty}^i}$ is a branched minimal immersion from $\Sigma_\infty^i=S^2$, hence $F_\infty$ is harmonic and conformal, and its Dirichlet energy can be estimated as 
$$
E(F_\infty) = \left| F_\infty(\S_\infty)\right|= \lim_{n\to\infty}  \left| F(\cdot,t_n)\circ\phi_n(S^2) \right| = \lim_{n\to\infty}\left|F(S^2,t_n)\right|\leq \left|F(S^2,0)\right|=E(w_k). 
$$
This contradicts the fact that there is no harmonic map from a stratified sphere to $M$ which is homotopic to $w_k$ but with energy not exceeding $E(w_k)$, by the choice of $w_k$.   \endproof

\section{Apendix}
In this section, we set $M=M_1\#M_2$, where $M_1$ and $M_2$ are 3-dimensional
closed manifold and $M$ is the connected sum of $M_1$ and $M_2$. That is to say,
we delete a small ball $B_1$ inside $M_1$ and a small ball $B_2$ inside $M_2$
and glue together their boundary spheres. We assume both $M_1$ and $M_2$ are
not $S^3$. 
We will use the following notations in this secton:
\begin{itemize}
\item[] $M_1'=M_1\backslash \overline{B_1}$,  $M_2'=M_2\backslash \overline{B_2}$, and
$S=\partial M_1'=\partial M_2'$.

\item[] $\widetilde{M}$: the universal cover of $M$ with covering map $\pi:\widetilde{M}\to M$.

\item[] $u$:  a diffeomorphism from $S^2$ to $S$.

\item[] $\tilde{u}$:  a lift of $u$.

\item[] $\gamma\circ \tilde{u}$: the deck transformation composed with $\tilde{u}$ for $\gamma\in\pi_1(M)$. 

\item[] For $v:S^2\rightarrow \widetilde{M}$,  $[v]_{H_2}$ denotes the element in $H_2(\widetilde{M})$ represented by $v$.

\end{itemize}

\begin{lem}\label{nontrivial}
 $u$ is non-trivial in $\pi_2(M)$.
\end{lem}

\proof Assume $u$ is trivial in $\pi_2(M)$. By Proposition 3.10 in \cite{Ha}, $u(S^2)$ bounds a compact contractible submanifold $B$ in $M$, and 
by the Poincar\'e conjecture, $B$ is a ball. Since $B\subset M\backslash S$, $S$
separates $M$ and $B$ is connected,  we have $B\subset M_1'$ or $B\subset M_2'$. We assume
$B\subset M_1'$. That  $M'_1\cap S=\emptyset$ implies 
$M_1'\backslash B=M_1'\backslash\overline{B}$ is open in $M'_1$. Since
$M_1'=(M_1'\backslash B)\cup B$ and $M_1'$ is connected, we conclude $B=M_1'$.
Thus $M_1 = B \cup \overline{B_1}=S^3$, which is
a contradiction to our assumption.
\endproof

\begin{lem}\label{no.intersection} For any $\gamma\in\pi_1(M)$ which is not 1, we have
$\gamma\circ \tilde{u}(S^2)\cap \tilde{u}(S^2)=\emptyset.$
\end{lem}

\proof Assume there is a point $p\in \gamma\circ \tilde{u}(S^2)\cap
\tilde{u}(S^2)$.
Let $\beta:[0,1]\rightarrow \widetilde{M}$ be a lift of $\gamma$ with $\beta(0)=p$.
Let $p_1=\gamma(1)\in\gamma\circ\tilde{u}(S^2)$ and $\beta_1:[0,1]\rightarrow
\gamma\circ\tilde{u}(S^2)$ be a continuous curve with $\beta_1(0)=p_1$ and $\beta_1(1)=p$.
Then the product path $\beta_2=\beta \cdot \beta_1$ is a loop in $\widetilde{M}$. Since $\pi_1(\widetilde{M})=\{1\}$,
 $\pi(\beta_2)$ is trivial in $\pi_1(M)$. The loop $\pi(\beta_1)$ is trivial in $\pi_1(M)$ as it is contained in $u(S^2)$. So 
$\pi(\beta)=1$  in $\pi_1(M)$, and it contradicts $\pi(\beta)= \gamma\not=1$.
\endproof

\begin{lem}\label{component}
 $\widetilde{M}\backslash \tilde{u}(S^2)$ has exactly
two connected components, and the closure of each of them is noncompact in
$\widetilde{M}$. Moreover, for any nontrivial $\gamma\in\pi_1(M)$,
$\widetilde{M}\backslash(\tilde{u}(S^2)\cup \gamma\circ\tilde{u}(S^2))$
has exactly 3 connected components, and only one of them has boundary $\tilde{u}(S^2)\cup\gamma\circ\tilde{u}(S^2)$.
\end{lem}

\proof
Since $\tilde{u}$ is an embedding, $\widetilde{M}\backslash \tilde{u}(S^2)$ is either connected or has two
components. Assume it is connected.
Let $\beta:[0,1]\rightarrow \widetilde{M}$ be a smooth curve in a neighborhood of $\tilde{u}(S^2)$ such that $\beta$ intersects 
$\tilde{u}$ at  $\beta(\frac{1}{2})$ transversally and $\beta(\frac{1}{2})$ is the only intersection
of $\tilde{u}$ and $\beta$. Since $\widetilde{M}\backslash \tilde{u}(S^2)$ is connected,
we can find a curve $\beta':[1,2]\rightarrow \widetilde{M}\backslash \tilde{u}(S^2)$
with $\beta'(1)=\beta(1)$ and $\beta'(2)=\beta(0)$. Let $\beta''(t)=\beta(t)$ for $t\in[0,1]$ and $\beta''(t)=\beta'(t)$ for $t\in[1,2]$. 
Then $\beta''$ is a loop which intersects $\tilde{u}(S^2)$ transversally at only one point.
Thus $\beta''$ is not trivial in $H_1(\widetilde{M})$, which contradicts the fact
that $\pi_1(\widetilde{M})=1$. So $\widetilde{M}\backslash\tilde{u}(S^2)$ has two components. 

Let $N_1$ and $N_2$ be two components of $\hat{M}\backslash \tilde{u}(S^2)$. Since
$\tilde{u}(S^2)\cap \gamma\circ\tilde{u}(S^2)=\emptyset$ by Lemma \ref{no.intersection}, we have $\gamma\circ\tilde{u}(S^2)\subset N_1$ or
$\gamma\circ\tilde{u}(S^2)\subset N_2$. Assume $\gamma\circ\tilde{u}(S^2)\subset N_2$. Using similar arguments, we
can prove $N_2\backslash \gamma\circ\tilde{u}(S^2)$ also has exact two components ($N_2$ is simply connected by the Van Kampen theorem). Obviously,
there is only one component whose boundary consists of  $\tilde{u}(S^2)$ and $\gamma\circ\tilde{u}(S^2)$.

If $\overline{N_1}$ is compact in $\widetilde{M}$, it is a smooth compact manifold with boundary $\tilde{u}(S^2)$, which implies $[\tilde{u}(S^2)]_{H_2(\widetilde{M})}=0$.
Hence $\tilde{u}(S^2)$ is trivial in $\pi_2(\widetilde{M})$ by the Hurewicz theorem.  However, it is well-known that $\pi_2(M)=\pi_2(\widetilde{M})$.
Lemma \ref{nontrivial} yields a contradiction. Similarly, $\overline{N_2}$ is noncompact in $\widetilde{M}$. 
\endproof

\begin{lem}
Each component $\hat{M}$ of $\pi^{-1}(M\backslash S)$ is a universal cover of $M_1'$ or $M_2'$.
For any $\gamma\in \pi_1(M)$, there are exactly two components, whose boundaries contain $\gamma\circ \tilde{u}(S^2)$.
\end{lem}

\proof Assume there is a point $p\in\hat{M}$, such that
$p'=\pi(p)\in M_1'$. We need to prove that $\hat{M}$ is a universal cover of $M_1'$.

First, we prove $\pi$ restricts to a covering map from $\hat{M}$ to $M_1'$. It suffices to
prove it is surjective. For any $p''\in M_1'$, we can find a curve $\beta:[0,1]
\rightarrow M_1'$ with $\beta(0)=p'$ and $\beta(1)=p''$. Let $\hat{\beta}$ be
the lift of $\beta$ in $\widetilde{M}$ with $\hat{\beta}(0)=p$. Since $\beta\subset M\backslash S$, $\hat{\beta}$ has no
intersection with $\pi^{-1}(S)$. Thus $\beta\subset \hat{M}$ as $\hat{M}$ is a connected component of $\pi^{-1}(M\backslash S)$. Hence $p''\in \pi(\hat{M})$.

Next, we show $\hat{M}$ is simple connected. Assume there is a close
curve $\delta:S^1\rightarrow \hat{M}$ which is not trivial in $\pi_1(\hat{M})$.
Then $\beta'=\pi(\delta)$
is a closed curve in $M_1'$. If 
$\beta'$ is nontrivial in $\pi_1(M_1')$, it is nontrivial in $\pi_1(M_1)$ as $M_1=M_1'\cup \overline{B_1}$. 
Since $\pi_1(M)=\pi_1(M_1)*\pi_1(M_2)$,  then $\beta'$ is also nontrivial in $\pi_1(M)$.
Note that $\delta$ is the lift of $\beta'$, so 
$\delta$ can not be a closed curve in $\hat{M}$. We get a contradiction.
If $\beta'$ is trivial in $\pi_1(M_1')$, since $\pi:\hat{M}\rightarrow M_1'$ is a covering map, $\delta$
must be trivial in $\hat{M}$. Again, we have a contradiction.
\endproof

\begin{lem}\label{H2}
Assume $M_2\neq \R\P^3$ and  $\hat{M}$ is a connected component of
$\pi^{-1}(M\backslash S)$, which is a universal cover of $M_2'$.  If $\gamma_1\neq
\gamma_2\in\pi_1(M)$, such that
$\gamma_1\circ \tilde{u}(S^2)$ and $\gamma_2\circ\tilde{u}(S^2)$ are in the
boundary of $\hat{M}$, then
$[\gamma_1\circ\tilde{u}]_{H_2}$ and
$[\gamma_2\circ\tilde{u}]_{H_2}$ are linearly independent in $H_2(\widetilde{M})$.
\end{lem}

\proof Since $M_2\not=S^3, \R\P^3$, we have $\pi_1(M_2)\neq \{1\}$ and
$\pi_1(M_2)\neq \mathbb{Z}_2$. Thus $\pi_1(M_2)$ contains at
least 3 elements. Then we can find $\gamma_3\in
\pi_1(M)$ which is not $\gamma_1,\gamma_2$, such that the closure of ${\hat{M}}$ contains at least 3 boundary components 
$\gamma_1\circ\tilde{u}(S^2)$, $\gamma_2\circ\tilde{u}(S^2)$ and
$\gamma_3\circ \tilde{u}(S^2)$.

Assume $[\gamma_1\circ \tilde{u}]_{H_2}$ and $[\gamma_2\circ\tilde{u}]_{H_2}$ are
linearly dependent. Equip $\widetilde{M}$ with the metric pulled back from $M$ via $\pi$. Then fixing a point $p\in \hat{M}$, we can find a sufficiently large $R$, such that the metric ball $B_{\frac{R}{2}}(p)$
in $\widetilde{M}$ contains $\gamma_i\circ\tilde{u}(S^2)$ for $i=1,2,3$, and
$[\gamma_1\circ\tilde{u}]_{H_2(B_R(p)\cap \widetilde{M})}$ and $[\gamma_2\circ\tilde{u}]_{H^2
(B_R(p)\cap \widetilde{M})}$ are linearly dependent in $H_2(B_R(p)\cap \widetilde{M})$.
We can choose $R$ via Sard's theorem and  by Lemma \ref{component},
such that the closure of ${B_R(p)\cap \widetilde{M}}$
is a smooth manifold with nonempty smooth boundary. 

By Lemma \ref{component},
$\widetilde{M}\backslash\{\gamma_1\circ\tilde{u}(S^2)\}$ has two components. Let  $N_1$
be the one which does not contain $\gamma_2\circ \tilde{u}(S^2)$ and
$N_3$ be the component of $\hat{M}\backslash
\gamma_3\circ\tilde{u}(S^2)$ which does not contain $\gamma_2\circ \tilde{u}(S^2)$.

Take a smooth curve $\beta:[0,1]\rightarrow {\hat{M}}\cup\partial\hat{M}$
so that $\beta$ intersects $\gamma_1\circ\tilde{u}(S^2)$ and
$\gamma_3\circ\tilde{u}(S^2)$ at $\beta(0)$ and $\beta(1)$ transversally, 
respectively, and demand $\beta((0,1))\cap \gamma_i\circ\tilde{u}(S^2)=
\emptyset$, $ i=1,2,3$. 
By Lemma \ref{component}, $N_1$ is unbounded, so $\partial B_R(p)\cap N_1\neq\emptyset$.
We can extend $\beta$ in $N_1$ to a curve
$\beta:[0,2]\rightarrow \widetilde{M}$ such that  $\gamma(2)\in\partial B_R$. Note that
$\gamma_2\circ\tilde{u}(S^2)\cap N_1=\emptyset$,
$\beta([0,2])\cap \gamma_2\circ\tilde{u}(S^2)=\emptyset$. Similarly, we can
extend $\beta$ to a curve $\beta:[-1,2]\rightarrow \widetilde{M}$, such that
$\beta(-1)\in\partial B_R$ and $\gamma([-1,2])\cap \gamma_2\circ\tilde{u}(S^2)
=\emptyset$.

Let $M^*$ be the manifold obtained by identifying two copies of $B_R(p)\cap \widetilde{M}$ along the boundary. Then $[\gamma_1\circ\tilde{u}]_{H_2(M^*)}$ and $[\gamma_2\circ\tilde{u}]_{H_2
(M^*)}$ are linearly dependent in $H_2(M^*)$. The two copies of $\beta$, however, yields a closed curve in $M^*$ which intersects $\gamma_1\circ\tilde{u}(S^2)$  transversally only at one point and has no intersection with $\gamma_2\circ\tilde{u}(S^2)$. This is impossible if $[\gamma_1\circ\hat{u}(S^2)]_{H_2}$ and $[\gamma_2\circ\hat{M}(S^2)]_{H_2}$ are linearly dependent. 
\endproof

\begin{lem}\label{cylinder}
If $[\gamma\circ\tilde{u}]_{H^2}=[\tilde{u}]_{H_2}$ where $\gamma\not=1\in\pi_1(M)$, then
$\widetilde{M}\backslash (\tilde{u}(S^2)\cup \gamma\circ\tilde{u}(S^2))$ has a component $\hat{M}(\gamma)$ which is topologically $S^2\times (-1,1)$ with
$\partial \hat{M}(\gamma)=\tilde{u}(S^2)\cup \gamma\circ\tilde{u}(S^2).$ 
\end{lem}

\proof
By Lemma \ref{no.intersection}, $\widetilde{M}\backslash (\tilde{u}(S^2)\cup
\gamma\circ\tilde{u}(S^2))$ has a component $\hat{M}(\gamma)$,
whose boundary consists of $\tilde{u}(S^2)$ and
$\gamma\circ\tilde{u}(S^2)$. Let $\hat{M}_1$ and $\hat{M}_2$
be the other two components with
$\partial \hat{M}_1=\tilde{u}(S^2)$ and $\partial\hat{M}_2=\gamma\circ\tilde{u}(S^2)$. 
We now prove the closure of ${\hat{M}(\gamma)}$ is compact.
Let $\beta:[-1,1]\rightarrow \widetilde{M}$ be a smooth curve which intersects
$\tilde{u}(S^2)$ at $\beta(0)$ transversally with
$\beta((0,1))\subset \hat{M}(\gamma)$ and $\beta((-1,0))\subset \hat{M}_1$. 
If the closure of ${\hat{M}(\gamma)}$ is not compact,  then fixing a point $p\in\hat{M}$,
for any large $R$,
$\partial B_R\cap \hat{M}(\gamma)\neq\emptyset$.
Extend $\beta|_{[0,1]}$ to a curve from $[0,2]$ to $\hat{M}(\gamma)\cap
\overline{B_R(p)}$ so that
$\beta(2)\in\partial B_R(p)\cap\hat{M}(\gamma),\beta((0,2))\subset \hat{M}(\gamma)$ and $ \beta((0,2))\cap \gamma\circ\tilde{u}(S^2)=\emptyset$. 
By Lemma \ref{component}, $\hat{M}_1$ is unbounded, then $\hat{M}_1\cap
\partial B_R(p)\neq\emptyset$. Similarly we extend $\beta|_{[-1,0]}$
in $\hat{M}_1$ to a curve $[-2,0]\rightarrow\hat{M}_1$ such that
$\beta(-2)\in\partial B_R(p)\cap\hat{M}_1, \beta([-2,0))\subset \hat{M}_1$, 
 which implies that
$\beta([-2,0))\cap \gamma\circ\tilde{u}(S^2)=\emptyset$. 
Then similar to the proof of Lemma \ref{H2}, we can construct a manifold
$M^*$, such that $[\tilde{u}(S^2)]_{H_2(M^*)}=[\gamma\circ\tilde{u}(S^2)]_{H_2(M^*)}$, and
there exists a closed curve $\beta'$, such that $\beta'$ intersect $\tilde{u}(S^2)$
in $M^*$ transversally at only one point and has no intersection with $\gamma\circ
\tilde{u}(S^2)$, which yields a contradiction. 

Next, we show $\hat{M}(\gamma)$ is simple connected.
Let $\hat{M}_i'$ be the manifold obtained by gluing a 3-ball to ${\hat{M}_i}$ along the boundary and let $\hat{M}'(\gamma)$ be the manifold
obtained by gluing a 3-ball to $\hat{M}(\gamma)$ along each boundary.
Since $\widetilde{M}=\hat{M}_1'\#\hat{M}'(\gamma)\#\hat{M}'_2$, 
$$\{1\}=\pi_1(\widetilde{M})=\pi_1(\hat{M}_1)*\pi_1(\hat{M}'(\gamma))
*\pi_1(\hat{M}'_2).$$
It follows $\pi_1(\hat{M}'(\gamma))=\{1\}$.
By the Poincar\'e conjecture, the closed simply connected $\hat{M}'(\gamma)$ is $S^3$ and in turn 
$\hat{M}(\gamma)=S^2\times (-1,1)$.
\endproof

\begin{lem}\label{not.homotopy} Let $M_2\neq \R \P^3$. If there exists $\gamma\not=1\in \pi_1(M)$ such that
$[\tilde{u}]_{H_2}=[\gamma\circ\tilde{u}]_{H_2}$,
then  $M_1=\R\P^3$ and $\gamma=\gamma'*1$, where $\gamma'$ is the nontrivial
element in $\pi_1(\R\P^3)$.
\end{lem}

\proof We claim that there is no  $\gamma_0\in\pi_1(M)$ that satisfiest $\gamma_0\neq 1$ or
$\gamma$ and
$\gamma_0\circ \tilde{u}(S^2)\subset \hat{M}(\gamma)$. Assume there is such a $\gamma_0$.
Let $\hat{M}$ be the component of
$\pi^{-1}(M\backslash S)$ which is a universal cover of $M_2'$
and $\gamma_0\circ\tilde{u}(S^2)$ belongs to $\partial\hat{M}$.
Lemma \ref{H2} asserts that $[\gamma\circ\tilde{u}(S^2)]_{H_2}$ and $[\gamma_0\circ\tilde{u}(S^2)]_{H_2}$ are linearly independent in $H_2(\hat{M}(\gamma))$.
But it follows from Lemma \ref{cylinder} that $H_2(\hat{M}(\gamma))=H_2(S^2\times (-1,1))=\mathbb{Z}$.
A contradiction.

Thus,
$\hat{M}(\gamma)$ is a connected component of $\pi^{-1}(M\backslash S)$. Thus, it must be
a universal cover of $M_1'$. Since $\hat{M}(\gamma)$ has only 2 boundary components,
we see $\pi_1(M_1)=\pi_1(M_1')=\mathbb{Z}_2$.
\endproof

\begin{pro}\label{h.sequence}
If $M_2\neq\R\P^3$, then we can find a sequence $\{\gamma_i:i=0,1,\cdots\}\subset \pi_1(M)$,
such that $[\gamma_i\circ u]_{H_2}\neq [\gamma_j\circ u]_{H_2}$ for any $i\neq j$.
\end{pro}

\proof Note that $\pi_1(M)=\pi_1(M_1)*\pi_1(M_2)$ contains infinitely many
elements since $\pi_1(M_1),\pi_1(M_2)\not=\{1\}$. By Lemma \ref{not.homotopy}, when $M_1\neq \R\P^3$,
for any sequence $\{\gamma_k\}\in\pi_1(M)$ with $\gamma_i\neq\gamma_j$ for any $i\neq j$, we
have $[\gamma_i\circ u]_{H_2}\neq [\gamma_j\circ u]_{H_2}$. 
When $M_1=\R\P^3$,  choose a sequence such that
$\gamma_i\neq \gamma_j\cdot\gamma'*1$ for any $i\neq j$, where $\gamma'$
is the nontrivial element of $\pi_1(\R\P^3)$.
By Lemma \ref{not.homotopy} again, $[\gamma_i\circ\tilde{u}]_{H_2}\neq
[\gamma_j\circ\tilde{u}]$ for any $i\neq j$.
\endproof

\end{document}